\documentclass{article}
\usepackage{hyperref}
\usepackage{arxiv}

\usepackage[utf8]{inputenc} 
\usepackage[T1]{fontenc}    
\usepackage{mathrsfs,amsmath,amssymb,bm}
\usepackage{lipsum}
\usepackage{amsfonts}
\usepackage{epstopdf}
\usepackage{makecell, rotating}
\usepackage{multirow}
\usepackage[nottoc]{tocbibind}
\usepackage{array}
\usepackage{listings}
\usepackage{cleveref}
\usepackage{caption}
\usepackage{authblk}
\usepackage[utf8]{inputenc} 
\usepackage[T1]{fontenc}    
\usepackage{url}            
\usepackage{booktabs}       
\usepackage{nicefrac}       
\usepackage{microtype}      
\usepackage{xcolor}         
\usepackage{multirow}
\usepackage{graphicx,subfigure}
\usepackage[ruled,linesnumbered]{algorithm2e}
\usepackage{tabularx}

\numberwithin{equation}{section}

\title{Fourier Neural Solver for large sparse linear algebraic systems}

\author[1]{Chen Cui}
\author[1]{Kai Jiang\thanks{\texttt{kaijiang@xtu.edu.cn}}}
\author[1]{Yun Liu}
\author[1]{Shi Shu\thanks{\texttt{shushi@xtu.edu.cn}}}

\affil[1]{Hunan Key Laboratory for Computation and Simulation in Science and Engineering, Key Laboratory of Intelligent Computing and Information Processing of Ministry of Education, School of Mathematics and Computational Science, Xiangtan University, Xiangtan, Hunan, China, 411105.}

\renewcommand{\headeright}{}
\renewcommand{\undertitle}{}

\begin{document}
\maketitle

\begin{abstract}
  Large sparse linear algebraic systems can be found in a variety of scientific and engineering fields, and many scientists strive to solve them in an efficient and robust manner. In this paper, we propose an interpretable neural solver, the Fourier Neural Solver (FNS), to address them. FNS is based on deep
  learning and Fast Fourier transform.  Because the error between the iterative solution and the ground truth involves a wide range of frequency modes, FNS combines a stationary iterative method and frequency space correction to eliminate different components of the error.  
  Local Fourier analysis reveals that the FNS can pick up on the error components in frequency space that are challenging to eliminate with stationary methods. Numerical experiments on the anisotropy diffusion equation, convection-diffusion equation, and Helmholtz equation show that FNS is more efficient and more robust than the state-of-the-art neural solver.
\end{abstract} 
	

\section{Introduction}\label{sec:01}
Large sparse linear algebraic systems are ubiquitous in scientific and engineering computation, such as discretization of partial differential equations (PDE) and linearization of nonlinear problems. Designing efficient, robust, and adaptive numerical methods for solving them is a long-term challenge. Iterative methods are an effective way to resolve this issue. They  can be classified into single-level and multi-level methods.
There are two types of single-level methods: stationary and non-stationary methods\,\cite{barrett1994templates}. 
Due to the sluggish convergence, stationary methods such as weighted Jacobi, Gauss-Seidel, successive over relaxation  methods\,\cite{saad2003iterative} are frequently utilized as smoothers in multi-level methods or as preconditioners.
Non-stationary methods typically refer to Krylov subspace methods, such as conjugate gradient (CG), generalized minimal residual (GMRES) methods\,\cite{hestenes1952methods,saad1986gmres}, whose convergence rate will be influenced heavily by some factors like initial value.
Multi-level methods mainly include geometric multigrid (GMG) method\,\cite{brandt1977multi,briggs2000multigrid, trottenberg2000multigrid} and algebraic multigrid (AMG) method\,\cite{falgout2006introduction,xu2017algebraic}. 
They are both composed of many factors, such as smoother and coarse grid correction, which heavily affect convergence. Finding these factors for a concrete problem is an art that demands extensive analysis, innovation, and trial.

In recent years, the technique of automatically picking parameters for Krylov and multi-level methods or constructing a learnable iterative scheme based on deep learning has attracted a lot of interest. For second-order elliptic equations with smooth coefficients, many neural solvers have achieved satisfactory results.
Hsieh et al.\,\cite{hsieh2019learning} utilized a convolutional neural network (CNN) to accelerate the convergence of Jacobi method. Luna et al.\,\cite{luna2021accelerating} accelerated the convergence of GMRES with a learned initial value. Significant efforts are also made in the development of multigrid solvers, such as learning smoother, transfer operator\,\cite{weymouth2021data,tomasi2021construction} and coarse-fine splitting\,\cite{taghibakhshi2021optimization}. 
For anisotropy elliptic equations, Huang et al.\,\cite{huang2021learning} exploited a CNN to design a more sensible smoother. Results showed that the magnitude of the learned smoother is dispersed along the anisotropic direction. Wang et al.\,\cite{wang4110904learning} introduced a learning-based local weighted least square method for the AMG interpolation operator, and applied it to random diffusion equations and one-dimensional small wavenumber Helmholtz equations.
Fanaskov\,\cite{fanaskov2021neural} learned the smoother and transfer operator of GMG in a neural network form. 
When the anisotropic strength is mild (within two orders of magnitude), previously mentioned works exhibit a considerable acceleration.
Chen et al.\,\cite{chen2022meta} proposed the Meta-MgNet to learn a basis vector of Krylov subspace as the smoother of GMG for strong anisotropic cases. However, the convergence rate is still sensitive to the anisotropic strength.
For convection-diffusion equations, Katrutsa et al.\,\cite{katrutsa2020black} learned the weighted Jacobi smoother and transfer operator of GMG, which has a  positive effect on the upwind discretazation system and also applied to solve a one-dimensional Helmholtz equation.
For second-order elliptic equations with random diffusion coefficients, Greenfeld et al.\,\cite{greenfeld2019learning} employed a residual network to construct the prolongation operator of AMG for uniform grids. Luz et al.\,\cite{luz2020learning} extended it to non-uniform grids using graph neural networks, which outperforms
classical AMG methods. For jumping coefficient problems, Antonietti et al.\,\cite{antonietti2021accelerating} presented a neural network to forecast the strong connection parameter to speed up AMG and used it as a preconditioner for CG. For Helmholtz equation, Stanziola et al.\,\cite{stanziola2021helmholtz} constructed a fully learnable neural solver, the helmnet, which is built on U-net and recurrent neural network\,\cite{kapturowski2018recurrent}. Azulay et al.\,\cite{azulay2022multigrid} developed a preconditioner based on U-net and shift-Laplacian MG\,\cite{erlangga2006novel} and applied the flexible GMRES\,\cite{calandra2012flexible} to solve the discrete system. For solid and fluid mechanics equations, there are also some neural solvers on associated discrete systems, such as but not limited to, learning initial values\,\cite{um2020solver,nikolopoulos2022ai}, constructing preconditioners\,\cite{stanaityte2020ilu}, learning search directions of CG\,\cite{kaneda2022deep}, learning parameters of GMG\,\cite{margenberg2022neural, margenberg2021deep}.

In this paper, we propose the Fourier Neural Solver (FNS), a deep learning and Fast Fourier Transform (FFT)\,\cite{cooley1965algorithm} based neural solver. FNS is made up of two modules: the stationary method and the frequency space correction.
Since stationary methods like weighted Jacobi method are difficult to get rid of low-frequency error, FNS uses FFT and CNN to learn these modes in the frequency space. 
Local Fourier analysis (LFA)\,\cite{brandt1977multi} reveals that FNS can pick up on the error components in frequency space that are challenging to eradicate by stationary methods. Therefore, FNS builds a complementary relationship by stationary method and CNN to eliminate error. With the help of FFT, the single-step iteration of FNS has a $O(N\log_2 N)$ computational complexity. All matrix-vector products are implemented using convolution, which is both storage-efficient and straightforward to parallelize. We investigate the effectiveness and robustness of FNS on three types of convection-diffusion-reaction equations.
For anisotropic equations, numerical experiments show that FNS can reduce the number of iterations by nearly $10$ times compared to the state-of-the-art Meta-MgNet when the anisotropic strength is relatively strong. For the non-symmetric systems arising from the convection-diffusion equations discretized by central difference method, FNS can converge while MG and CG methods diverge. And FNS is faster than other algorithms such GMRES and  BiCGSTAB($\ell$)\,\cite{sleijpen1993bicgstab}. For the indefinite systems arising from the Helmholtz equations, 
FNS outperforms GMRES and BiCGSTAB at medium wavenumbers.
In this paper, we apply FNS to above three PDE systems. However, the principles employed by FNS indicate that FNS has the potential to be useful for a broad range of sparse linear algebraic systems.

The rest of this paper is organized as follows. Section \ref{sec:02} proposes a general form of linear convection-diffusion-reaction equation and describes the motivation for designing FNS. Section \ref{sec:03} presents the FNS algorithm. Section \ref{sec:04} examines the performance of FNS to anisotropy, convection-diffusion, and Helmholtz equations. Finally, Section \ref{sec:05} draws the conclusions and future work.

\section{Motivation}\label{sec:02}
We consider the general linear convection-diffusion-reaction equation with Dirichlet boundary condition 
\begin{equation}
    \begin{cases}-\varepsilon \nabla \cdot(\boldsymbol{\alpha}(x) \nabla u)+\nabla \cdot(\boldsymbol{\beta}(x) u)+ \gamma u =f(x) & \text { in } \Omega \\ u(x)=g(x) & \text { on } \partial \Omega\end{cases}
    \label{eq:general}
\end{equation}
where $\Omega \subseteq \mathbb{R}^d$ is an open and bounded domain. $\boldsymbol{\alpha}(x)$ is the $d\times d-$ order diffusion coefficient matrix. $\boldsymbol{\beta}(x)$ is the $d\times 1$ velocity field that the quantity is moving with. $\gamma$ is the reaction coefficient. $f$ is the source term.

We can obtain a linear algebraic system once we discretize Eq.\,\eqref{eq:gexiangyixing} by finite element method (FEM) or finite difference method (FDM)
\begin{equation}
     \mathbf{A} \mathbf{u}=\mathbf{f},
     \label{eq:linear}
\end{equation}
where $\mathbf{A} \in \mathbb{R}^{N \times N}$,  $\mathbf{f} \in \mathbb{R}^{N} $ and $N$ is the spatial discrete degrees of freedom.

Classical stationary iterative methods, such as Gauss-Seidel and weighted Jacobi methods, have the generic form
\begin{equation}
    \mathbf{u}^{k+1}=\mathbf{u}^{k}+\mathbf{B}\left(\mathbf{f}-\mathbf{A} u^{k}\right),
    \label{eq:simple_iter}
\end{equation}
where $\mathbf{B}$ is an easily computed operator such as the inverse of diagonal matrix (Jacobi method), the inverse of lower triangle matrix (Gauss-Seidel method).
However, the convergence rate of such methods is relatively low. As an example, we utilize the weighted Jacobi method to solve a special case of Eq.\,\eqref{eq:general} and use LFA to analyze the reason.

Taking $\varepsilon=1,\,\boldsymbol{\alpha}(x)=\left(\begin{array}{ll}1 & 0\\ 0 & 1\end{array}\right)$, $\boldsymbol{\beta}(x)=\left(\begin{array}{l}0\\ 0\end{array}\right)$ and $\gamma=0$, Eq.\,\eqref{eq:general} becomes the Poisson equation.
With a linear FEM discretization and in stencil notation, the resulting discrete operator reads
\begin{equation}
   \left[\begin{array}{ccc}
    & -1 & \\
    -1 & 4 & -1\\
    &-1 &
    \end{array}\right].
    \label{eq:poisson_stencil}
\end{equation}
In weighted Jacobi method, $\mathbf{B}=\omega\mathbf{I}/4$, where $\omega \in (0,1]$ and $\mathbf{I}$ is the identity matrix. Eq. \eqref{eq:simple_iter} can be written in the pointwise form 
\begin{equation}
u^{k+1}_{ij}=u^{k}_{ij}+\frac{\omega}{4}\left(f_{ij}-(4u^k_{ij}-u^{ k}_{i-1, j}-u^{k}_{i+1, j}-u^{k}_{i, j-1}-u^{k}_{i, j+ 1})\right).
\end{equation}
Let $u_{ij}$ be the true solution, and define error $e^k_{ij} = u_{ij} - u^{k}_{ij}$. Then we have 
\begin{equation}
    e^{k+1}_{ij}=e^{k}_{ij}-\frac{\omega}{4}(4e^k_{ij} - e^{k}_{i-1, j} - e^{k}_{i+1, j} - e^{k}_{i, j-1} - e^{k}_{i, j+1}).
    \label{eq:errorpoint}
\end{equation}

Expanding error  in a Fourier series $e^k_{ij}=\sum_{p_1,p_2} v^{k} e^{i 2 \pi\left(p_{1} x_{i}+p_{2} y_{j}\right)}$, substituting the general term $v^{k} e^{i 2 \pi\left(p_{1} x_{i}+p_{2} y_{j}\right)}$, $p_{1}, p_{2} \in[-N/2, N/2)$ into Eq.\,\eqref{eq:errorpoint}, we have
\begin{equation}\nonumber
    \begin{aligned}
        v^{k+1} &=v^{k}\left(1-\frac{\omega}{4}\left(4-e^{-i 2 \pi p_{1} h}-e^ {i 2 \pi p_{1} h}-e^{-i 2 \pi p_{2} h}-e^{-i 2 \pi p_{2} h}\right)\right)\\
    &= v^{k}\left(1-\frac{\omega}{4}\left(4- 2\cos \left(2 \pi p_{1} h\right)-2 \cos \left( 2 \pi p_{2} h\right)\right)\right).
\end{aligned}
\end{equation}
The convergence factor of weighted Jacobi method (also known as smoother factor in MG framework\,\cite{trottenberg2000multigrid}) is
\begin{equation}
    \mu_\text{loc}:= \left|\frac{v^{k+1}}{v^{k}}\right|=\left|1-\omega+\frac{\omega}{2}\left(\cos \left(2 \pi p_{1} h\right)+ \cos \left(2 \pi p_{2} h\right)\right)\right|.
\end{equation}
\begin{figure}[!htbp]
    \centering
    \subfigure[Distribution of convergence factor $\mu_\text{loc}$]{\label{subfig:poisson_factor}\includegraphics[width=0.4\textwidth]{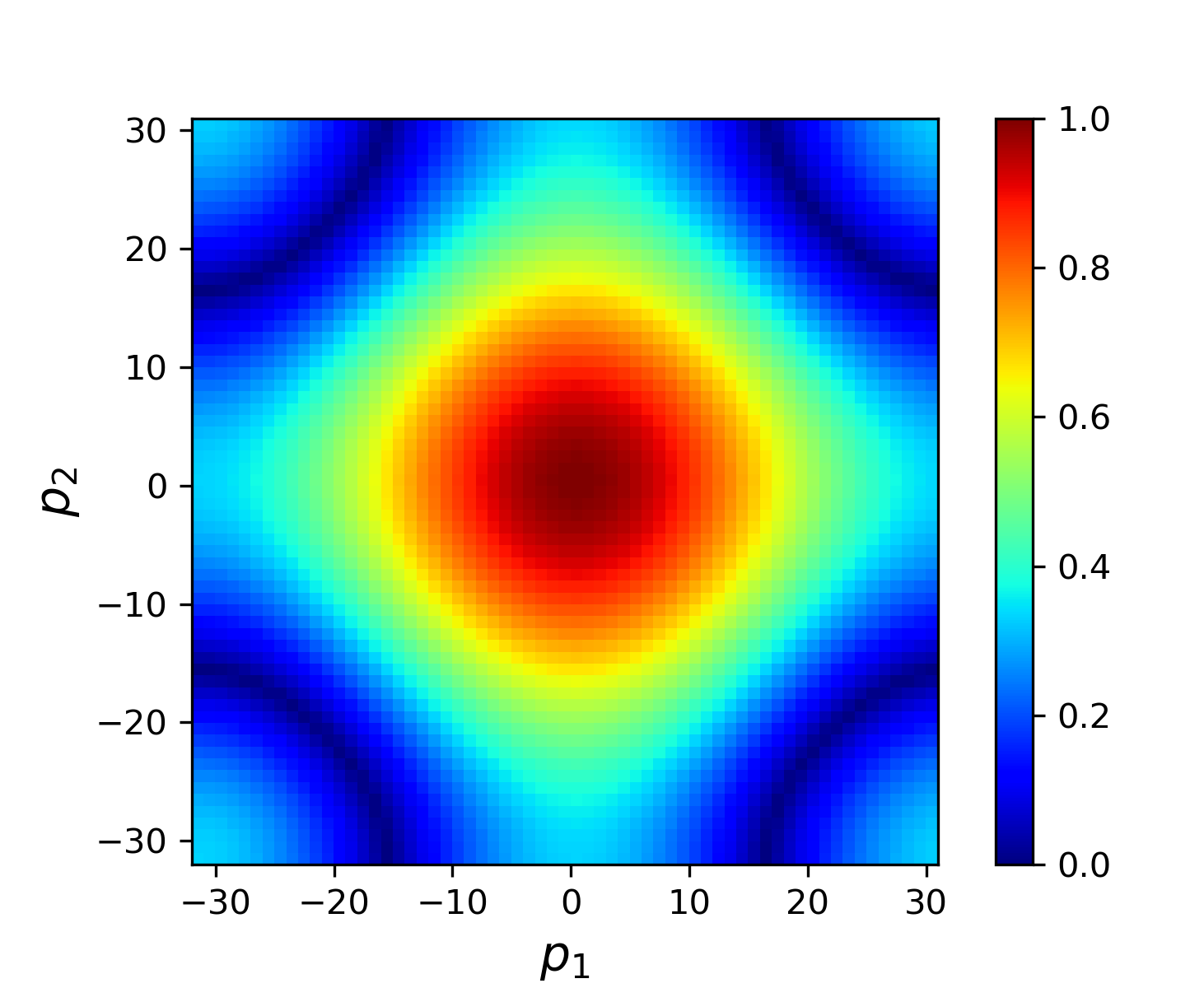}}\quad
    \subfigure[High and low-frequency regions]{\label{subfig:highlow}\includegraphics[width=0.4\textwidth]{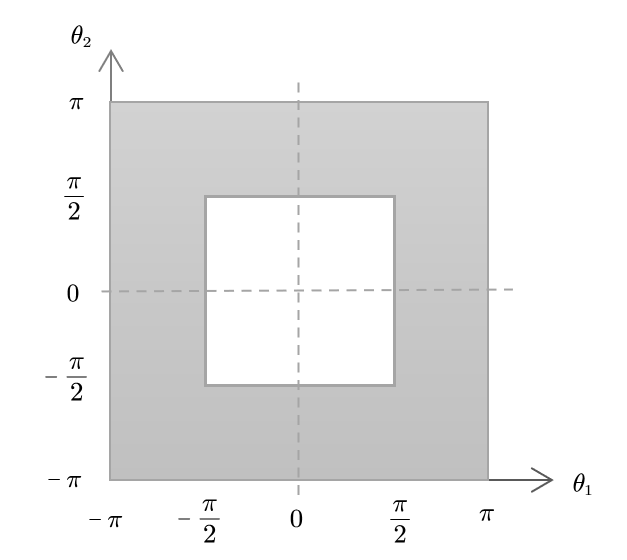}}
    \caption{(a) The distribution of convergence factor $\mu_\text{loc}$ of weighted Jacobi method ($\omega=2/3$) in solving linear system arised by Poisson equation; (b) Low-frequency (white) and high-frequency (gray) regions.}
\end{figure}

Figure~\ref{subfig:poisson_factor} shows the distribution of convergence factor $\mu_\text{loc}$ of weighted Jacobi ($\omega=2/3$) in solving linear system for Poisson equation.
For a better understanding, let $\theta_1 = 2\pi p_1 h$, $\theta_2 = 2\pi p_2 h$, $\boldsymbol{\theta} = (\theta_1, \theta_2) \in[ -\pi, \pi)^2$. Define the high and low-frequency regions 
\begin{equation}
    \begin{aligned}
    &T^{\text {low }}:=\left[-\frac{\pi}{2}, \frac{\pi}{2}\right)^{2}, \\
    & T^{\text {high }}:=[-\pi, \pi)^{2} \backslash\left[-\frac{\pi}{2}, \frac {\pi}{2}\right)^{2},
    \end{aligned}
\end{equation}
as shown in Figure~\ref{subfig:highlow}. It can be seen that in the high-frequency region, $\mu_\text{loc}$ is approximately zero, whereas in the low-frequency region, $\mu_\text{loc}$ is close to one. As a result, weighted Jacobi method attenuates high-frequency errors quickly but is mostly useless towards low-frequency errors. 

The reason is attributed to two aspects. Firstly, the high-frequency reflects the local oscillation, while the low-frequency reflects the global pattern. Since $\mathbf{A}$ is sparse and the basic operation $\mathbf{Au}$ of weighted Jacobi method is a local operation, which makes it challenging to remove low-frequency global error components. Secondly, $\mathbf{A}$ is sparse, $\mathbf{A}^{-1}$ is commonly dense, which means that the mapping $\mathbf{f}\rightarrow \mathbf{A}^{-1}\mathbf{f}$ is non-local, making local operations difficult to approximate.

Therefore, we should seek the solution in another space to obtain an effective approximation of the non-local mapping. For example, the Krylov methods approximate the solution in a subspace spanned by a set of basis. MG generates a coarsen space to broaden the receptive field of the local operation. 
However, as mentioned in Section \ref{sec:01}, these methods have various parameters to be carefully designed. In this paper, we propose the FNS, a generic solver that uses FFT to learn solutions in frequency space, with the parameters automatically obtained in a data-driven manner.

\section{Fourier Neural Solver}\label{sec:03}

Denote stationary iterative methods of 
\eqref{eq:simple_iter} in a operator form
\begin{equation}
    \mathbf{v}^{k+1}=\Phi(\mathbf{u}^{k}),
\end{equation}
and the $k-$th step residual
\begin{equation}
    \mathbf{r}^{k}:=\mathbf{f}-\mathbf{A}\mathbf{v}^{k+1},
\end{equation}
then the $k-$th step error $\mathbf{e}^k:=\mathbf{u}-\mathbf{v}^{k+1}$ satisfies residual equation
\begin{equation}
    \mathbf{A}\mathbf{e}^k=\mathbf{r}^{k}.
    \label{eq:residual}
\end{equation}

As shown in the preceding section, the slow convergence rate of stationary methods is due to the difficulty in reducing low-frequency errors. In fact, even high-frequency errors might be not effectively eliminated by $\Phi$ for many cases. Now, we employ stationary methods to rapidly erase some components of the error and use FFT to convert the remaining error components to frequency space. The resulting solver is the Fourier Neural Solver.

\begin{figure}[h]
  \centering
  \includegraphics[width=13.5cm]{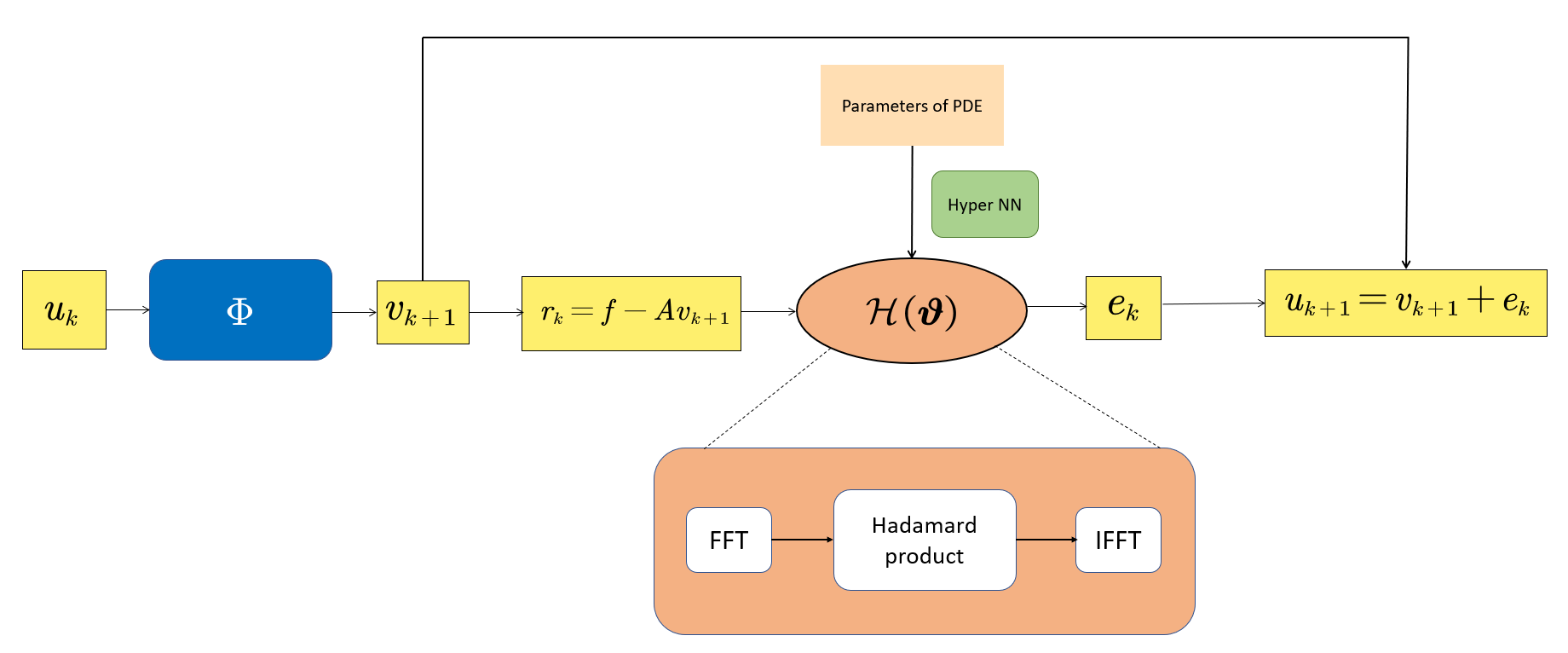}
  \caption{FNS calculation flow chart.}
  \label{fig:model}
\end{figure}

Figure~\ref{fig:model} shows a flowchart of the $k-$th step of FNS.
The module for solving the residual equation in frequency space is denoted as $\mathcal{H}$. It consists of three steps: FFT$\rightarrow$Hadamard product$\rightarrow$IFFT.
The parameter $\boldsymbol{\vartheta}$ of $\mathcal{H}$ is the output of the Hyper-neural network (HyperNN). The input $\boldsymbol{\eta}$ of HyperNN is the PDE parameters corresponding to the discrete systems. During training, the only parameter $\boldsymbol{\theta}$ of HyperNN serves as our optimization parameter.

The three-step operation of $\mathcal{H}$ is inspired by the fast Poisson solver\,\cite{swarztrauber1977methods}. Let eigenvalues and eigenvectors of $\mathbf{A}$ be $\lambda_{1}, \ldots, \lambda_{N}$ and $\mathbf{q}_{1}, \ldots, \mathbf{q}_{N}$, respectively. Solving Eq.\eqref{eq:linear} includes three steps:
\begin{enumerate}
    \item Expand $\mathbf{f}$ as a combination $\mathbf{f}=a_{1} \mathbf{q}_{1}+\cdots+a_{N} \mathbf{q}_{N} $ of the eigenvectors
    \item Divide each $a_{k}$ by $\lambda_{k}$
    \item Recombine eigenvectors into $\mathbf{u}=\left(a_{1} / \lambda_{1}\right) \mathbf{q}_{1}+\cdots+\left(a_{N} / \lambda_ {N}\right) \mathbf{q}_{N}$.
\end{enumerate}
In particular, when $\mathbf{A}$ is the system generated by five-point stencil \eqref{eq:poisson_stencil}, its eigenvector $\mathbf{q}_k$ is the sine function. The above first and third steps can now be done with a computational complexity of $O(N \log _{2} N)$ using Fast Sine Transform (based on the FFT). The computational complexity of each iteration of FNS is $O(N \log _{2} N)$. 

It is worth noting that, although $\Phi$ can smooth some components of the error, the components that are removed are indeterminate. As a result, instead of filtering high-frequency modes in frequency space, $\mathcal{H}$ learns error components that $\Phi$ cannot easily eliminate. For $\Phi$ with a fixed stencil, we can use LFA to demonstrate that the learned $\mathcal{H}$ is complementary to $\Phi$.

The loss function used here for training is the relative residual
\begin{equation}
     \mathcal{L} =\sum_{i=1}^{N_{b}} \frac{\|\mathbf{f}_i-\mathbf{A}_i\mathbf{u}_i^K\|_2} {\|\mathbf{f}_i\|_2},
     \label{eq:loss}
\end{equation}
where $\{\mathbf{A}_i,\mathbf{f}_i\}$ is the training data. $N_b$ is the batch size. $K$ indicates that the $K-$th step solution $\mathbf{u}^K$ is used to calculate the loss. These specific values will be given in the next section. The training and testing algorithms of FNS are summarized in Algorithm\,\ref{algorithm:FNS_offline} and Algorithm\,\ref{algorithm:FNS_online}, respectively.

\begin{algorithm}
    \caption{FNS offline traning}\label{algorithm:FNS_offline}
    \KwData{PDE parameters $\{\boldsymbol{\eta}_i\}_{i=1}^{N_{train}}$ and corresponding discrete systems $\{\mathbf{A}_i,\mathbf{f}_i\}_{i=1}^{N_{train}}$}
    \KwIn{$\Phi$, HyperNN($\boldsymbol{\theta}$), K and Epochs} 
    \For(serial){epoch $ = 1,\ldots,$ Epochs}{
    \For(parallel){$i=1,\ldots, N_{train}$}{
        $\boldsymbol{\vartheta}_i = $ HyperNN $(\boldsymbol{\eta}_i,  \boldsymbol{\theta})$

        $\mathbf{u}_i^0=$ zeros like $\mathbf{f}_i$

        \For(serial){$k=0,\ldots, K-1$}{
   
           $\mathbf{v}_i^{k+1}=\Phi(\mathbf{u}_i^{k})$

           $\mathbf{r}_i^{k}=\mathbf{f}_i-\mathbf{A}_i\mathbf{v}_i^{k+1}$

           $\mathbf{\hat{r}}_i^{k}=\mathcal{F}(\mathbf{r}_i^{k})$

           $\mathbf{\hat{e}}_i^{k}=\mathbf{\hat{r}}_i^{k}\circ \boldsymbol{\vartheta}_i$

           $\mathbf{e}_i^{k}=\mathcal{F}^{-1}(\mathbf{e}_i^{k})$

           $\mathbf{u}_i^{k+1}=\mathbf{v}_i^{k+1}+\mathbf{e}_i^{k}$
        }
    }
    Compute loss function \eqref{eq:loss}

    Apply Adam algorithm\,\,\cite{kingma2014adam} to update parameters $\boldsymbol{\theta}$
    }
    \Return{learned FNS}
\end{algorithm}

\begin{algorithm}[H]
    \caption{FNS online testing}\label{algorithm:FNS_online}
    \KwData{PDE parameter $\boldsymbol{\eta}$ and corresponding discrete system $\mathbf{A},\mathbf{f}$}
    \KwIn{Learned FNS, acceptable tolerance $tol$ and maximum number of iteration steps $MaxIterNum$} 
        Setup: $\boldsymbol{\vartheta}=$ HyperNN $(\boldsymbol{\eta}, \boldsymbol{\theta})$

        $k=0$

        $\mathbf{u}^0=$ zeros like $\mathbf{f}$

        $res = \frac{\|\mathbf{f}-\mathbf{A}\mathbf{u}^k\|_2}{\|\mathbf{f}\|_2}$

        \While{$res >$ $tol$ and $k <$ MaxIterNum}{
   
           $\mathbf{v}^{k+1}=\Phi(\mathbf{u}^{k})$

           $\mathbf{r}^{k}=\mathbf{f}-\mathbf{A}\mathbf{v}^{k+1}$

           $\mathbf{\hat{r}}^{k}=\mathcal{F}(\mathbf{r}^{k})$

           $\mathbf{\hat{e}}^{k}=\mathbf{\hat{r}}^{k}\circ \boldsymbol{\vartheta}$

           $\mathbf{e}^{k}=\mathcal{F}^{-1}(\mathbf{e}^{k})$

           $\mathbf{u}^{k+1}=\mathbf{v}^{k+1}+\mathbf{e}^{k}$

           $res = \frac{\|\mathbf{f}-\mathbf{A}\mathbf{u}^k\|_2}{\|\mathbf{f}\|_2}$

          $k=k+1$
        }
    
    \Return{solution of $\mathbf{A} \mathbf{u}=\mathbf{f}$}
\end{algorithm}

\section{Numerical Experiments}\label{sec:04}
We take the anisotropy equation, convection-diffusion equation, and Helmholtz equation as examples to demonstrate the performance of FNS. In all experiments, matrix-vector products are implemented by CNN based on Pytorch\,\cite{paszke2019pytorch} platform. All code can be found at \url{https://github.com/cuichen1996/FourierNeuralSolver}.

\subsection{Anisotropy equation}
Consider the anisotropy diffusion equation
\begin{equation}
    \left\{\begin{aligned}
    -\nabla \cdot(C \nabla u) &=f, \quad \text { in } \Omega , \\
    u &=0, \quad \text { on } \partial \Omega ,
    \end{aligned}\right.
    \label{eq:gexiangyixing}
\end{equation}
the diffusion coefficient matrix 
\begin{equation}
    C=C(\xi  , \theta)=\left(\begin{array}{cc}\cos \theta & -\sin \theta \\ \sin \theta & \cos \theta\end{array}\right)\left(\begin{array}{ll}1 & 0 \\ 0 & \xi  \end{array}\right)\left(\begin{array}{cc}\cos \theta & \sin \theta \\ -\sin \theta & \cos \theta\end{array}\right),
\end{equation}
$0<\xi <1$ is the anisotropic strength, $\theta \in[0, \pi]$ is the anisotropic direction, $\Omega=(0,1) ^2$.
We use bilinear FEM to discretize \eqref{eq:gexiangyixing} with a uniform $n\times n$ quadrilateral mesh. The associated discrete system is shown in \eqref{eq:linear}, where $N = (n- 1)\times (n-1)$. We will carry out experiments for the following two cases.

\subsubsection{Case 1: Generalization ability of anisotropic strength with fixed direction}\label{subsec:41}

In this case, we use the same training and testing data as\,\cite{chen2022meta}. For fixed $\theta = 0$, we randomly sample $20$ distinct parameters $\xi$ from the distribution $\log_{10} \frac{1}{\xi } \sim U[0,5]$ and obtain $\{\mathbf{A}_i\}_{i=1}^{20}$ by discretizing \eqref{eq:gexiangyixing} using bilinear FEM with $n=256$. We randomly select $100$ right-hand functions for each $\mathbf{A}_i$. And each entry of $\mathbf{f}$ is sampled from the Gaussian distribution $N(0, 1)$. Therefore there are $N_{train}=2000$ training data.  The hyperparameters used for training including batch size, learning rate, $K$ in the loss function, and the concrete network structure of HyperNN are listed in Appendix\,\ref{appen:ani}.

FNS can take various kinds of $\Phi$. In this case, we use weighted Jacobi, Chebyshev semi-iterative (Cheby-semi), and Krylov methods. The weight of weighted Jacobi method is $2/3$. Krylov method uses the same subspace as\,\cite{chen2022meta}. Its basis vector is the output of a DenseNet\,\cite{huang2017densely}. For Chebyshev semi-iterative method, we provide a brief summary here. More details can refer to\,\cite{golub1961chebyshev, adams2003parallel}.

\begin{figure}[!htbp]
    \centering
    \subfigure{\includegraphics[width=0.18\textwidth]{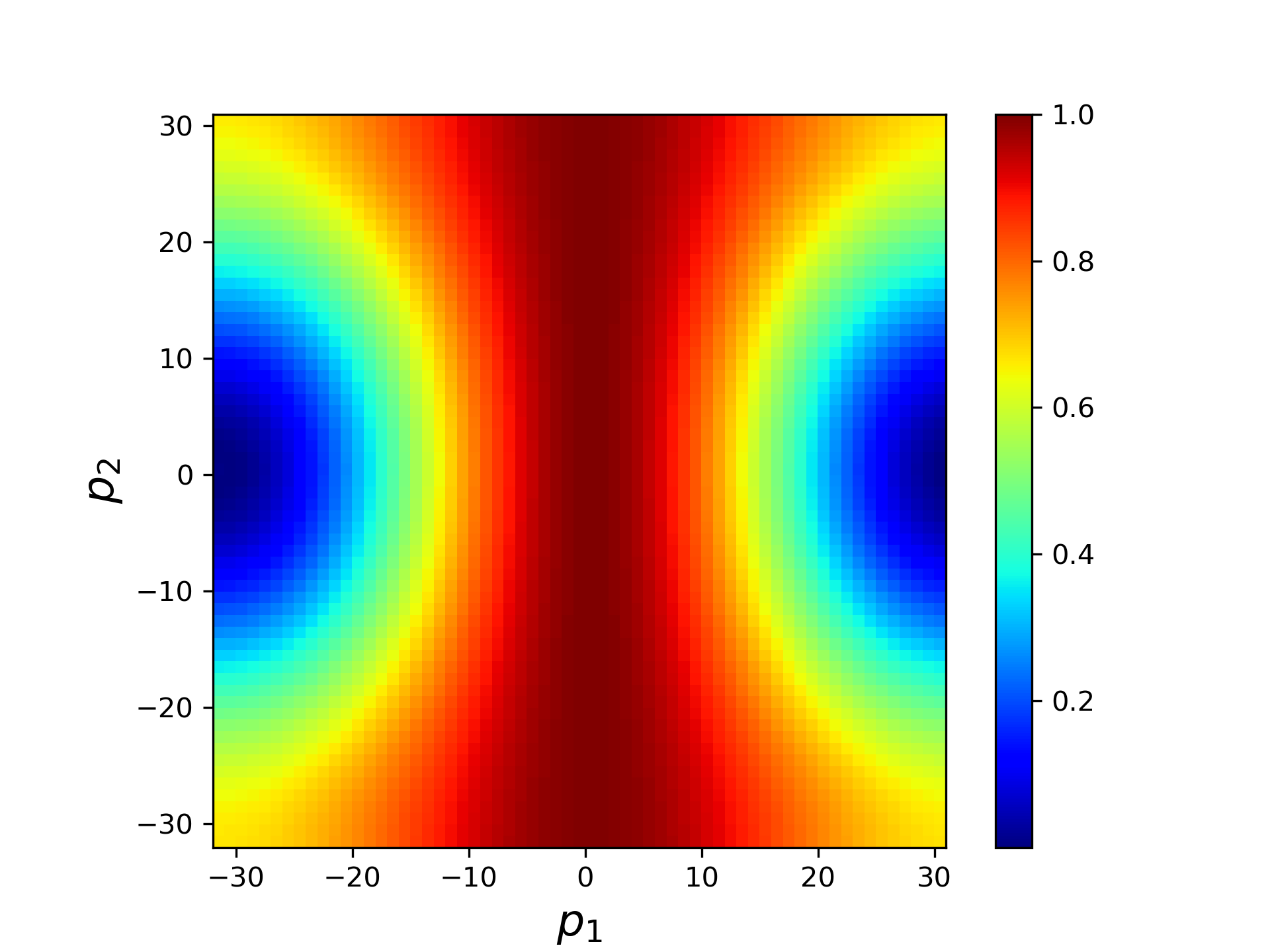}}
    \subfigure{\includegraphics[width=0.18\textwidth]{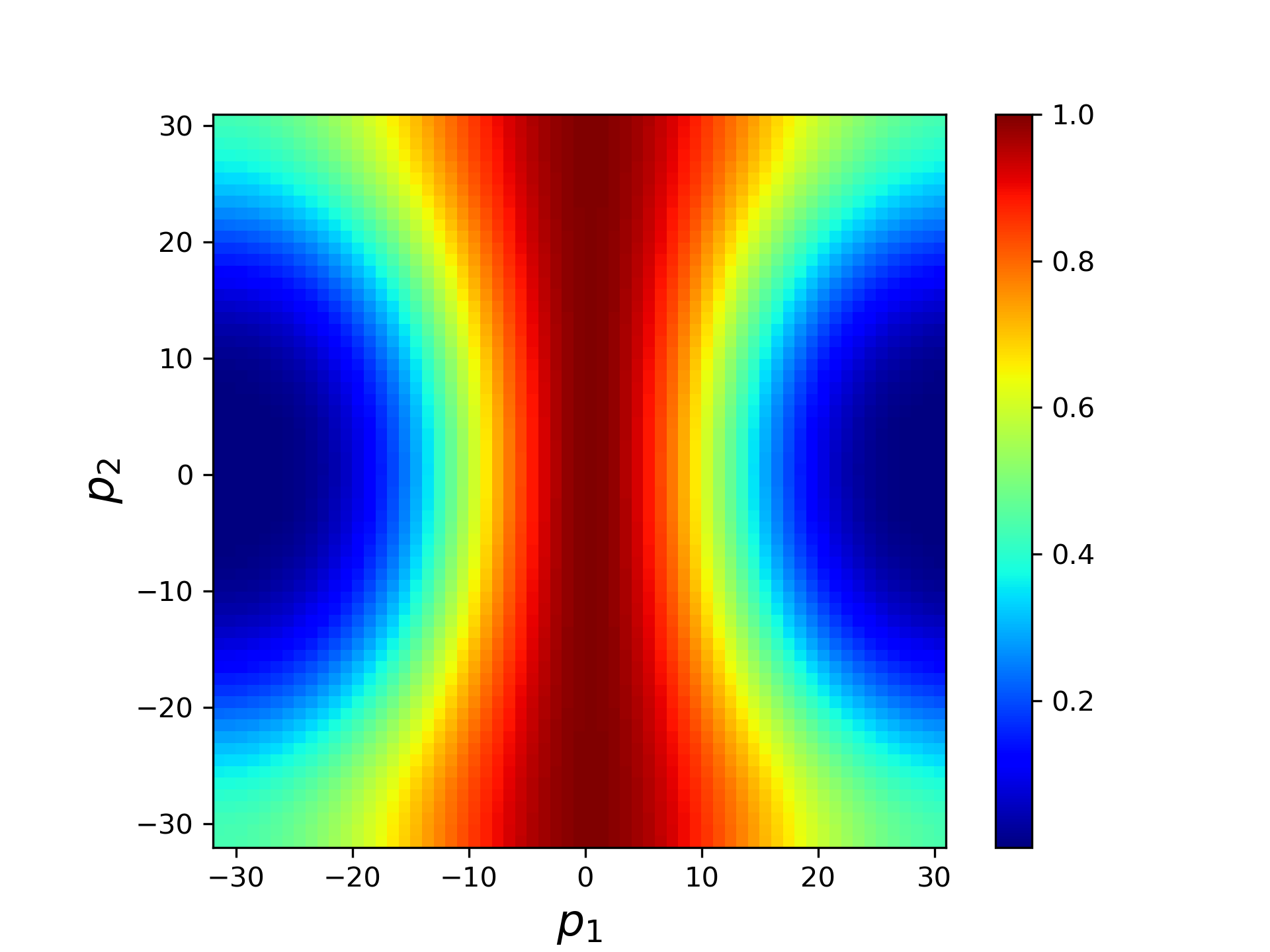}}
    \subfigure{\includegraphics[width=0.18\textwidth]{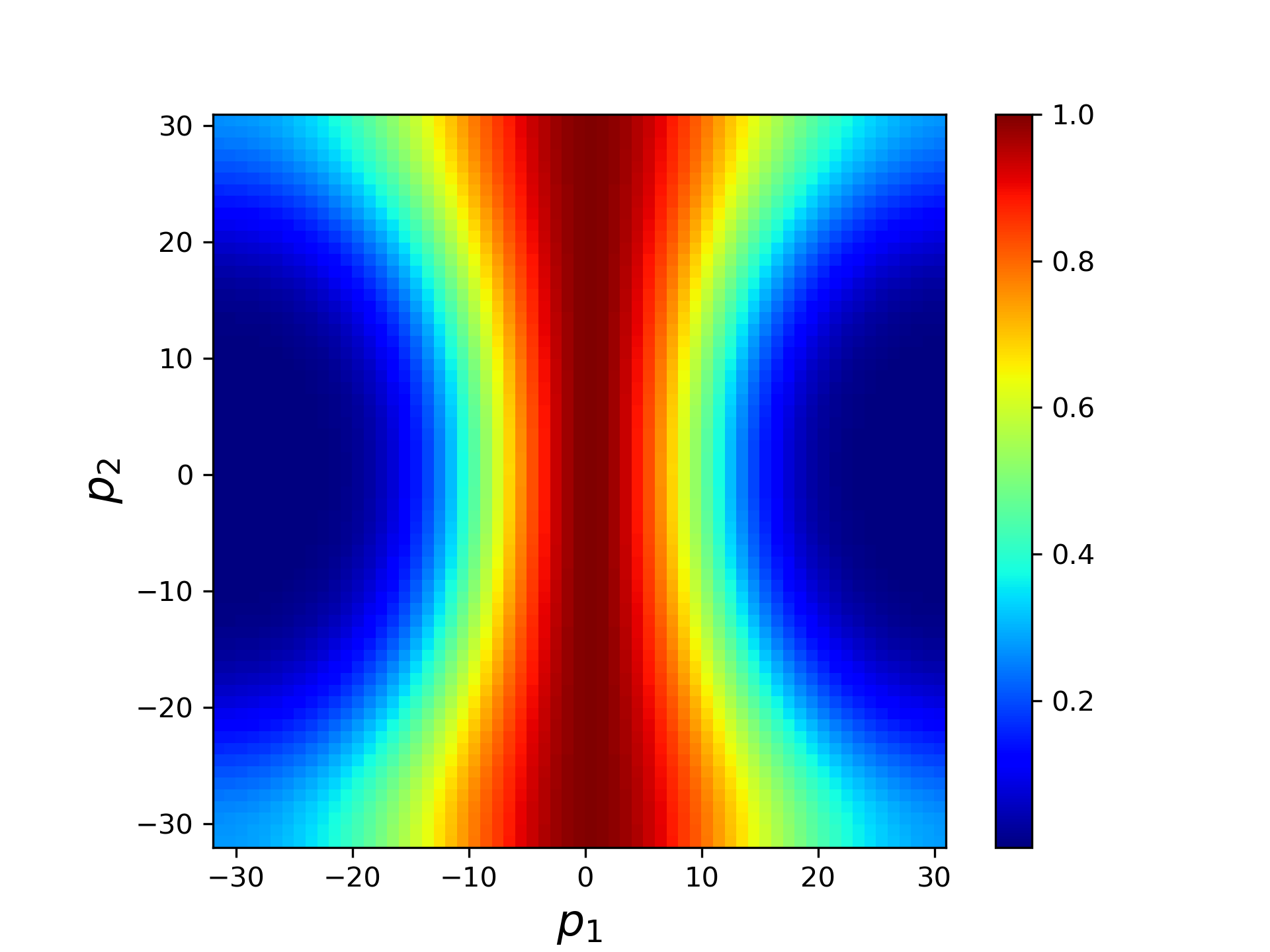}}
    \subfigure{\includegraphics[width=0.18\textwidth]{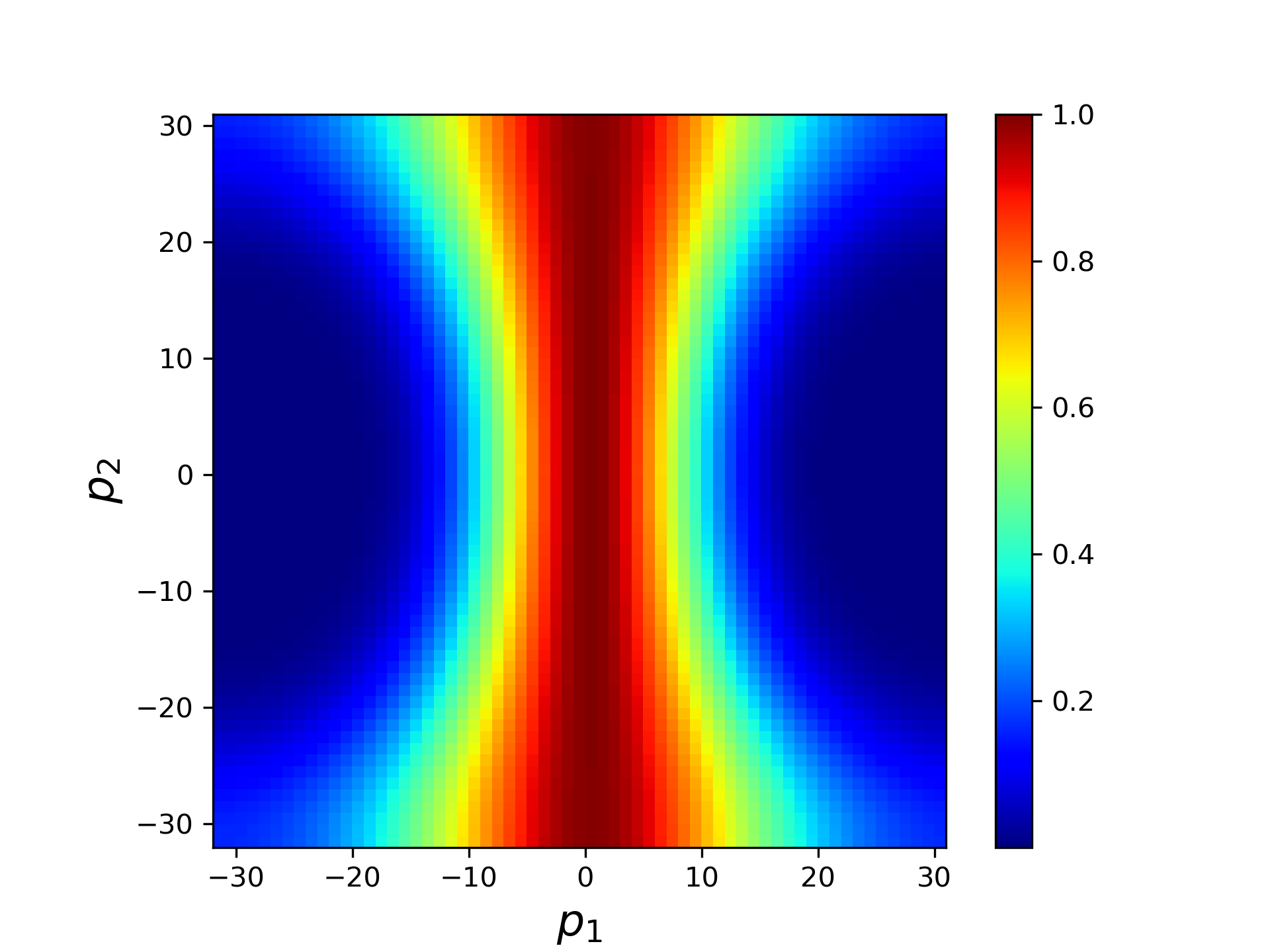}}
    \subfigure{\includegraphics[width=0.18\textwidth]{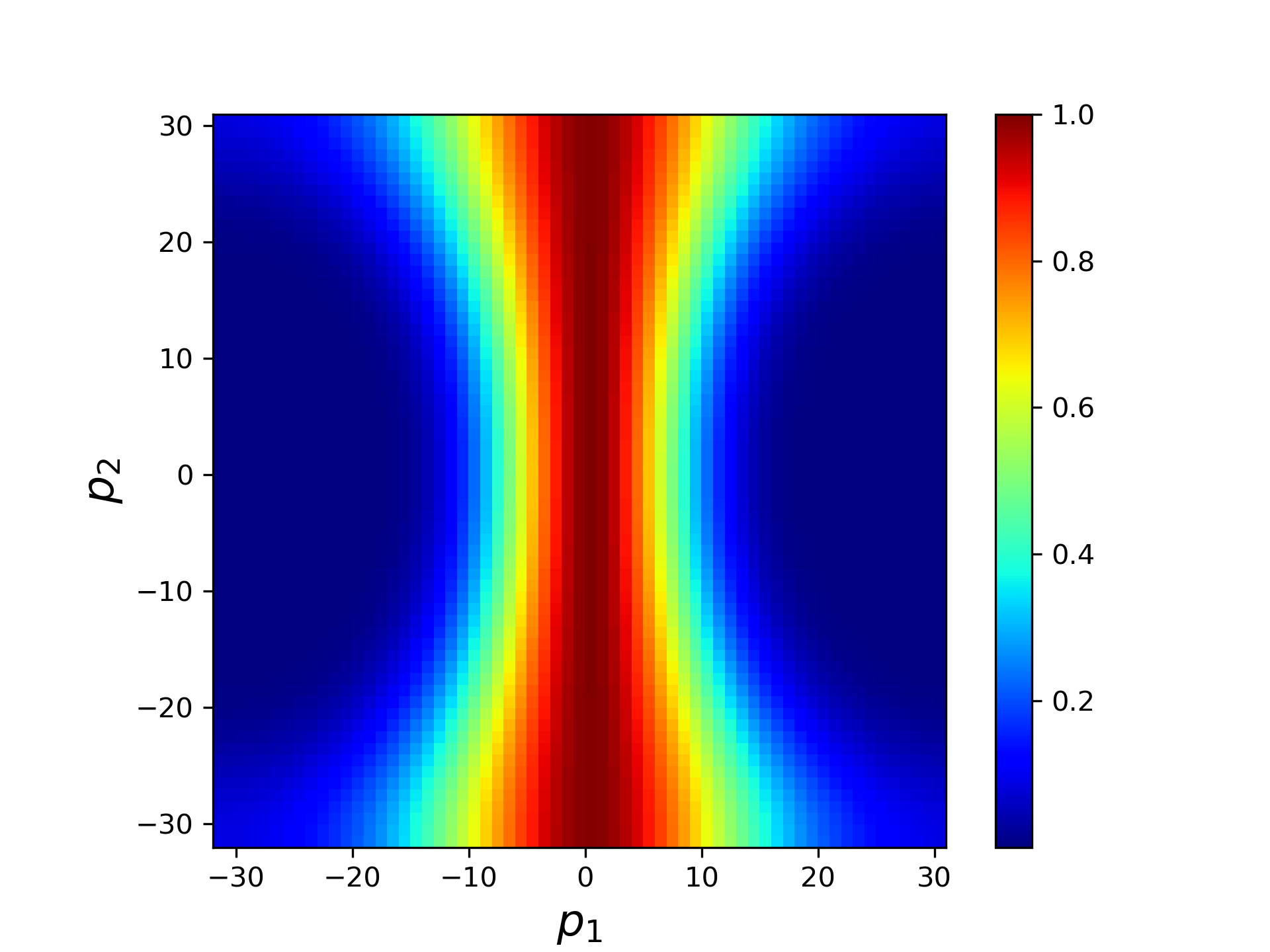}}\\
    \subfigure{\includegraphics[width=0.18\textwidth]{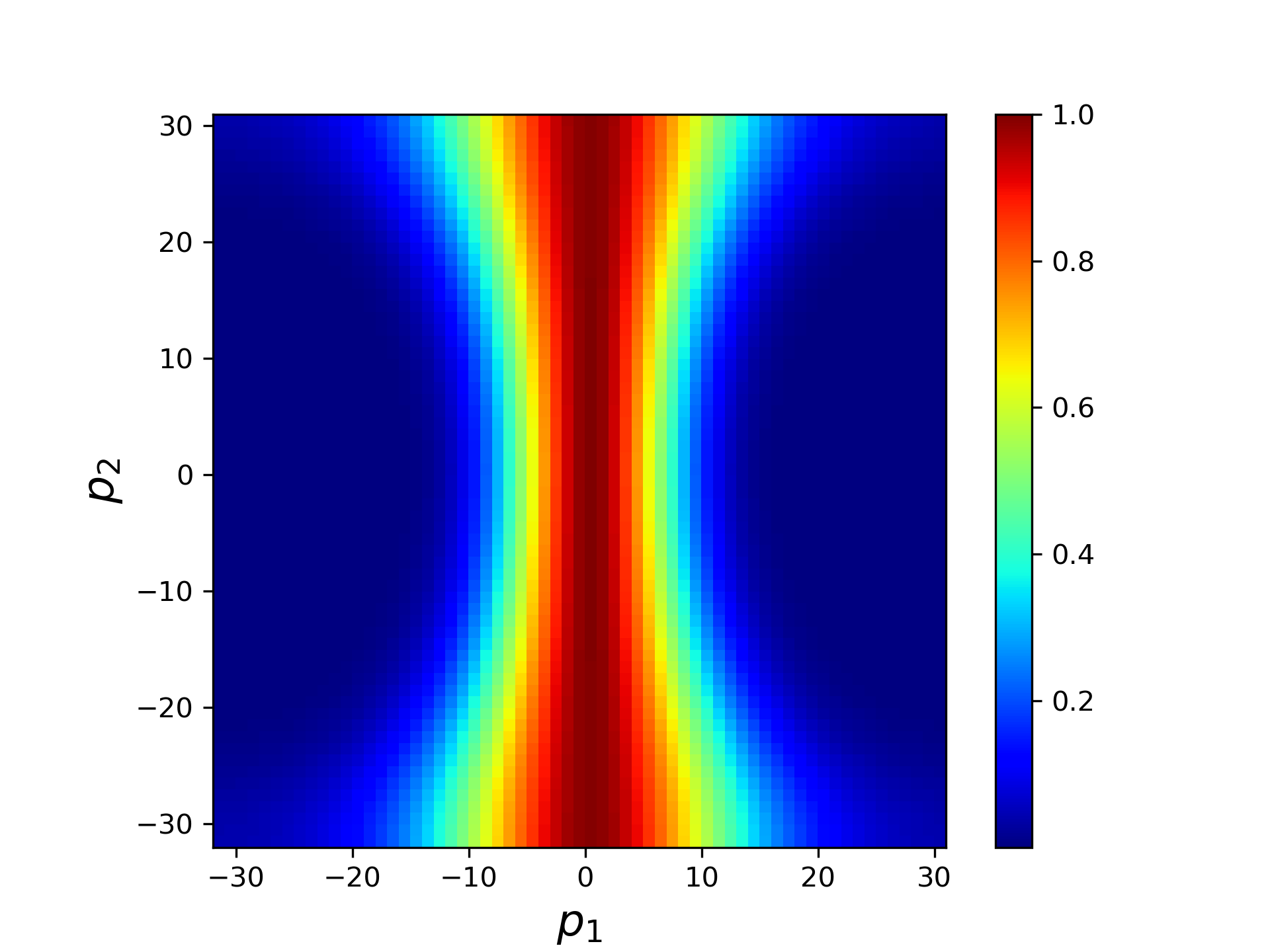}}
    \subfigure{\includegraphics[width=0.18\textwidth]{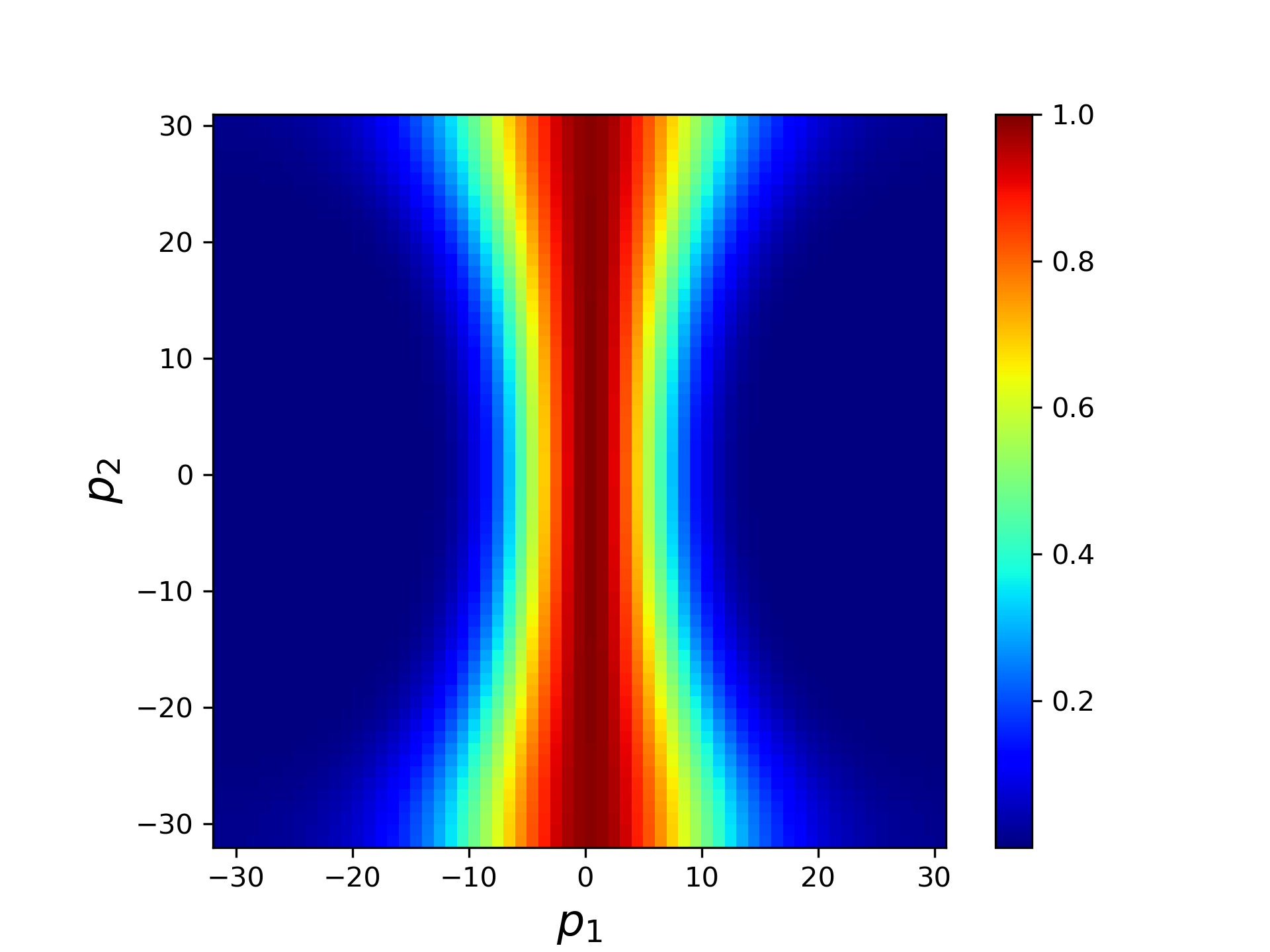}}
    \subfigure{\includegraphics[width=0.18\textwidth]{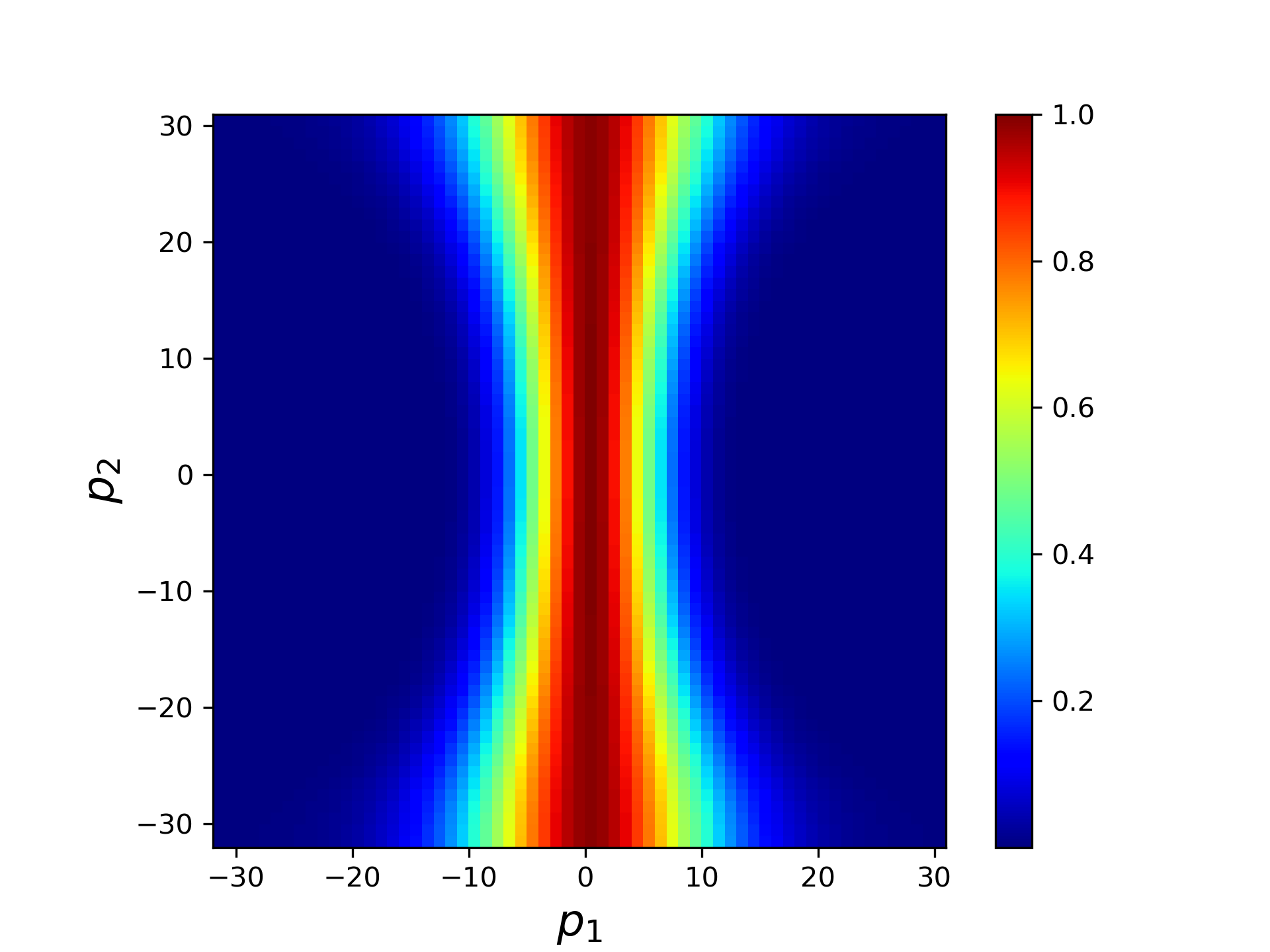}}
    \subfigure{\includegraphics[width=0.18\textwidth]{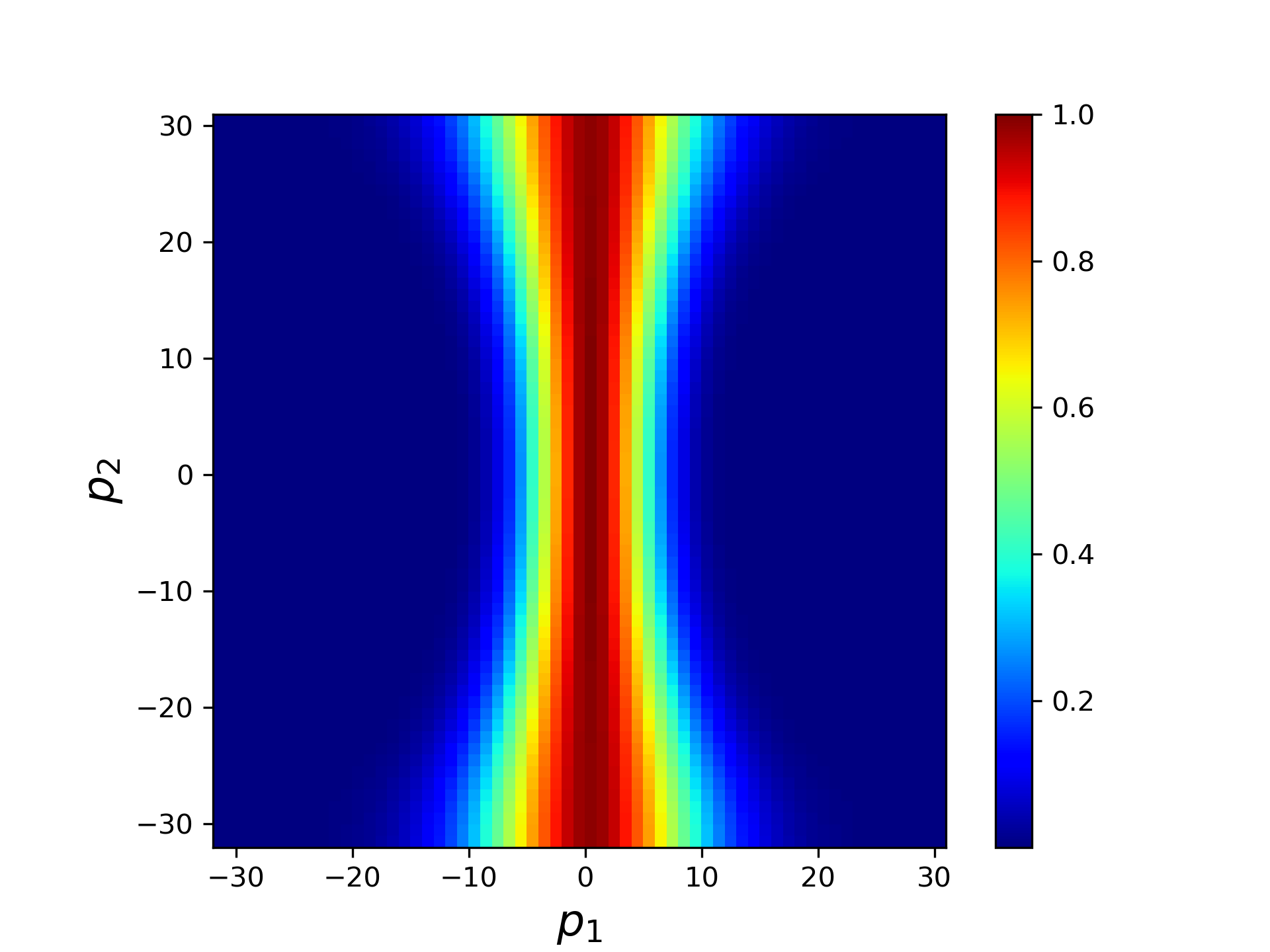}}
    \subfigure{\includegraphics[width=0.18\textwidth]{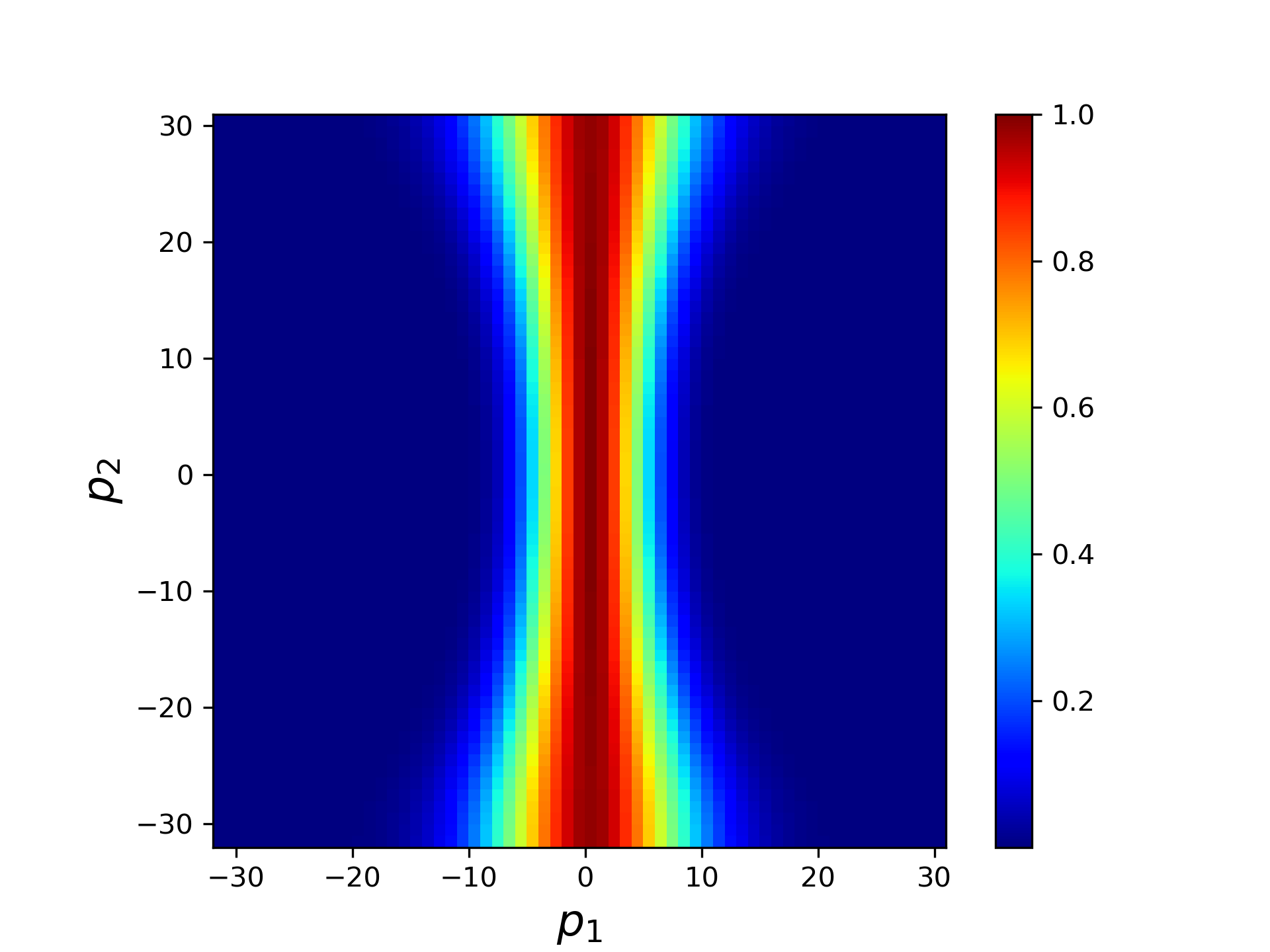}}
    \caption{The distribution of convergence factor for Cheby-semi ($m=10$) when $\xi=10^{-3}, \theta=0$.}
    \label{fig:ChebyLFA}
\end{figure}

If we vary the parameter of Richardson iteration at each step
\begin{equation}
    \mathbf{u}^{k+1}=\mathbf{u}^{k}+\tau_k\left(\mathbf{f}-\mathbf{A} u^{k}\right),
\end{equation}
and the maximum and minimum eigenvalues of $\mathbf{A}$ are known, then $\tau_k$ can be determined as
\begin{equation}
\tau_{k}=\frac{2}{\lambda_{\max }+\lambda_{\min }-\left(\lambda_{\min }-\lambda_{\max }\right) x_{k}} , \quad k=0, \ldots, m-1,
\end{equation}
where
\begin{equation}
x_{k}=\cos \frac{\pi(2 k+1)}{2 n}, k=0, \ldots, m-1,
\end{equation}
are the roots of a $m-$order Chebyshev polynomial.
Here, $\lambda_{\max}$ is obtained by the power method\,\cite{mises1929praktische}, but calculating $\lambda_{\min }$ often incurs an expensive computational cost. Therefore, we use $\lambda_{\max}/\alpha$ to replace $\lambda_{\min }$. The resulting method  is referred to as Chebyshev semi-iteration. We take $m=10, \alpha=3$. Figure~\ref{fig:ChebyLFA} shows the convergence factor obtained by LFA. It can be observed that the smooth effect improves as $m$ increases. However, the high-frequency error along $y$ direction is also difficult to eliminate. When $\Phi$ is Jacobi method, we implement $\Phi$ five-times, and transform residual to frequency space to correct error. This is because employing the stationary method multi-times can enhance its smoothing effect.

After training, we choose $\theta=0,\,\xi = 10^{-1},\ldots, 10^{-6}$ for testing and generate $10$ random right-hand functions for each parameter. The iteration is terminated when the relative residual is less than $10^{-6}$. We use "mean $\pm$ std" to show the mean and standard deviation of iterations over test set as \cite{chen2022meta} does.

Table~\ref{tab:test_param_case1} shows test results of different solvers. It can be found that iteration steps of all solvers grow as anisotropic strength increases except for MG (line-Jacobi). The growth of FNS is substantially lower than Meta-MgNet and MG. When FNS employs the same $\Phi$ as Meta-MgNet, the number of iterations is nearly $10$ times lower than that of Meta-MgNet at $\varepsilon=10^{-5}$ with the same computational complexity of single-step iteration. It is also worth mentioning that the line-Jacobi smoother can only be applied to several specific $\theta$, i.e. $0, \pi/4, \pi/2, 3\pi/4, \pi$. However, FNS can be available for arbitrary $\theta$.
\begin{table}[h]
  \caption{The mean and standard deviation of the number of iterations required to achieve the stopping criterion over all tests for the anisotropy equation case 1. "$-$" means that it cannot converge within 10000 steps, and "$\,\,$" means that \,\cite{chen2022meta} does not provide test results for this parameter.}
  \label{tab:test_param_case1}
  \resizebox{\textwidth}{!}{
  \begin{tabular}{|c|cccccc|}
      \hline $\xi$ & FNS(Cheby-semi) & FNS(Jacobi) & FNS(Krylov)&  Meta-MgNet(Krylov)\cite{chen2022meta} & MG(Jacobi)& MG(line-Jacobi)  \\
      \hline
      $\xi =10^{-1}$ & $67.9\pm 3.81$ & $138.9\pm 11.18$ & $30.0\pm 4.58$ & $7.5 \pm 0.50$ & $90.2\pm 0.98$ & $13.0 \pm 0.00$\\
      $\xi =10^{-2}$ & $101.6 \pm 8.72$ & $167.8\pm 13.81$ & $38.5\pm 3.83$ & $35.1\pm 1.04$ & $752.8 \pm 12.23$  & $13.0 \pm 0.00 $\\
      $\xi =10^{-3}$ & $151.0\pm 7.24$ &  $221.7\pm 11.56$& $48.6\pm 3.26$ & $171.6 \pm 6.34$ & $5600\pm119.42$ & $13.0 \pm 0.00$  \\
      $\xi =10^{-4}$ & $233.2 \pm 5.67$ & $330.1\pm 9.16 $& $65.5 \pm 2.80$ & $375.2 \pm 5.88$ & $-$  & $11.0 \pm 0.00$\\
      $\xi =10^{-5}$ & $340.1 \pm 9.43$ &$466.2\pm 13.47$ & $80.7\pm 7.21$ & $797.8 \pm 12.76$ & $-$ & $11.0 \pm 0.00 $\\
      $\xi =10^{-6}$ & $348.1 \pm 11.15$ & $477.9\pm 16.10$& $85.9 \pm 7.52$ & &$-$ &\\  
      \hline
  \end{tabular}}
\end{table}

\begin{figure}[!htbp]
    \centering
    \subfigure[$\xi = 10^{-1}$: convergence factor of $\Phi$]{
        \label{subfig:ACase1_2_a}
        \includegraphics[width=0.27\textwidth]{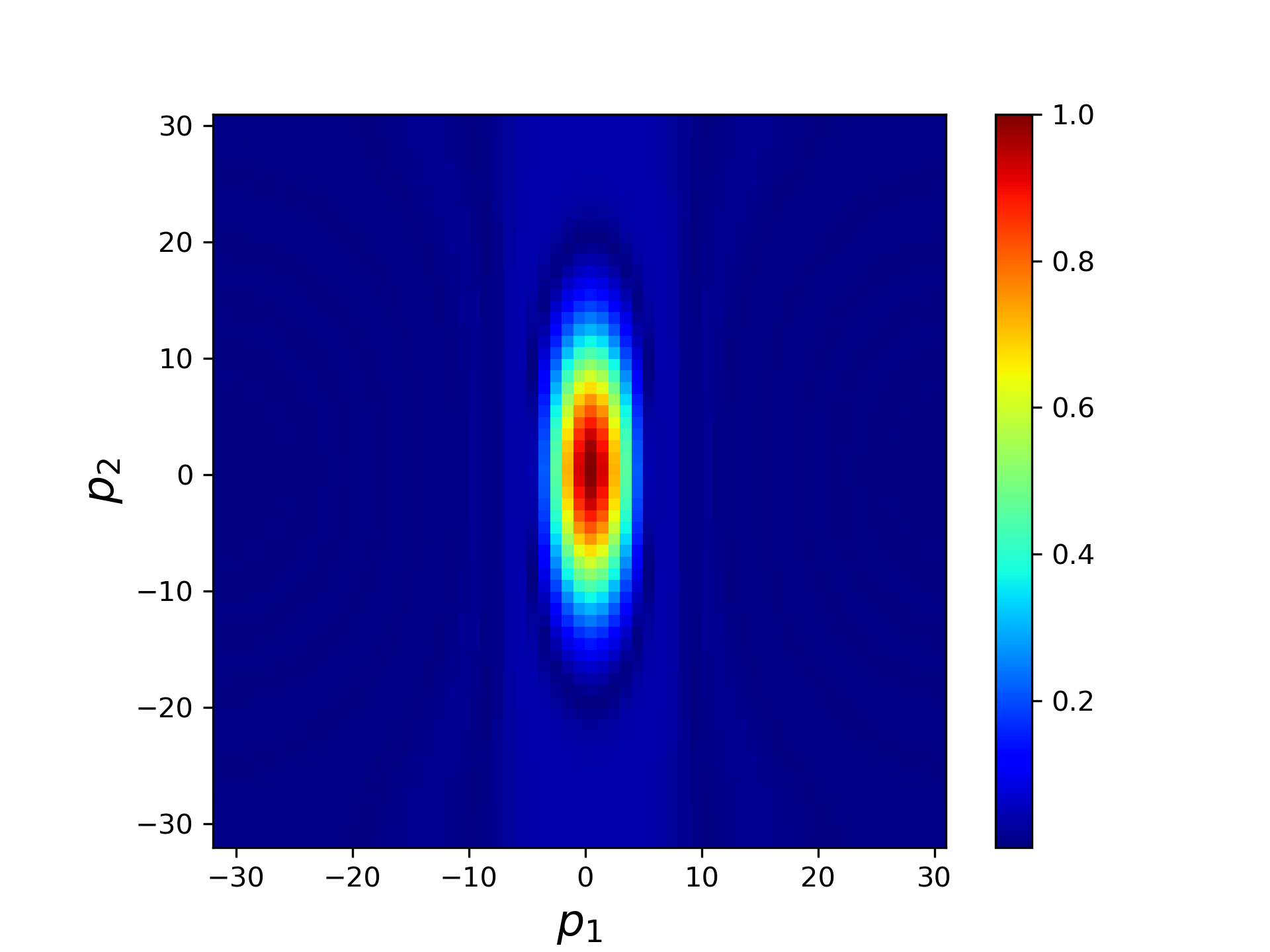}
    }
    \subfigure[$\xi = 10^{-1}$: $\hat{\mathbf{e}}$, before doing corrections]{
            \label{subfig:ACase1_2_b}
            \includegraphics[width=0.27\textwidth]{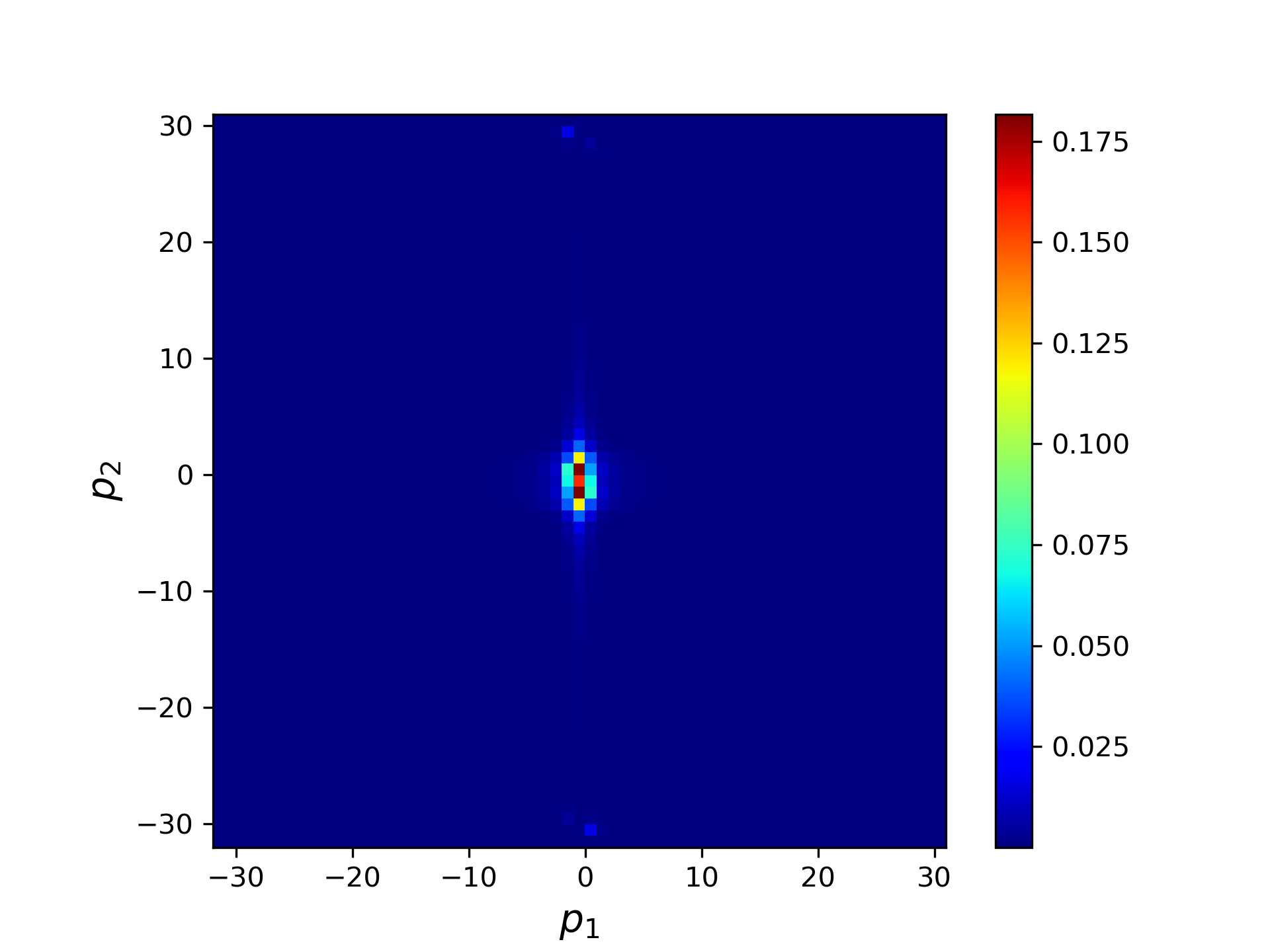}
    }
    \subfigure[$\xi = 10^{-1}$: learned $\hat{\mathbf{e}}$]{
            \label{subfig:ACase1_2_c}
            \includegraphics[width=0.27\textwidth]{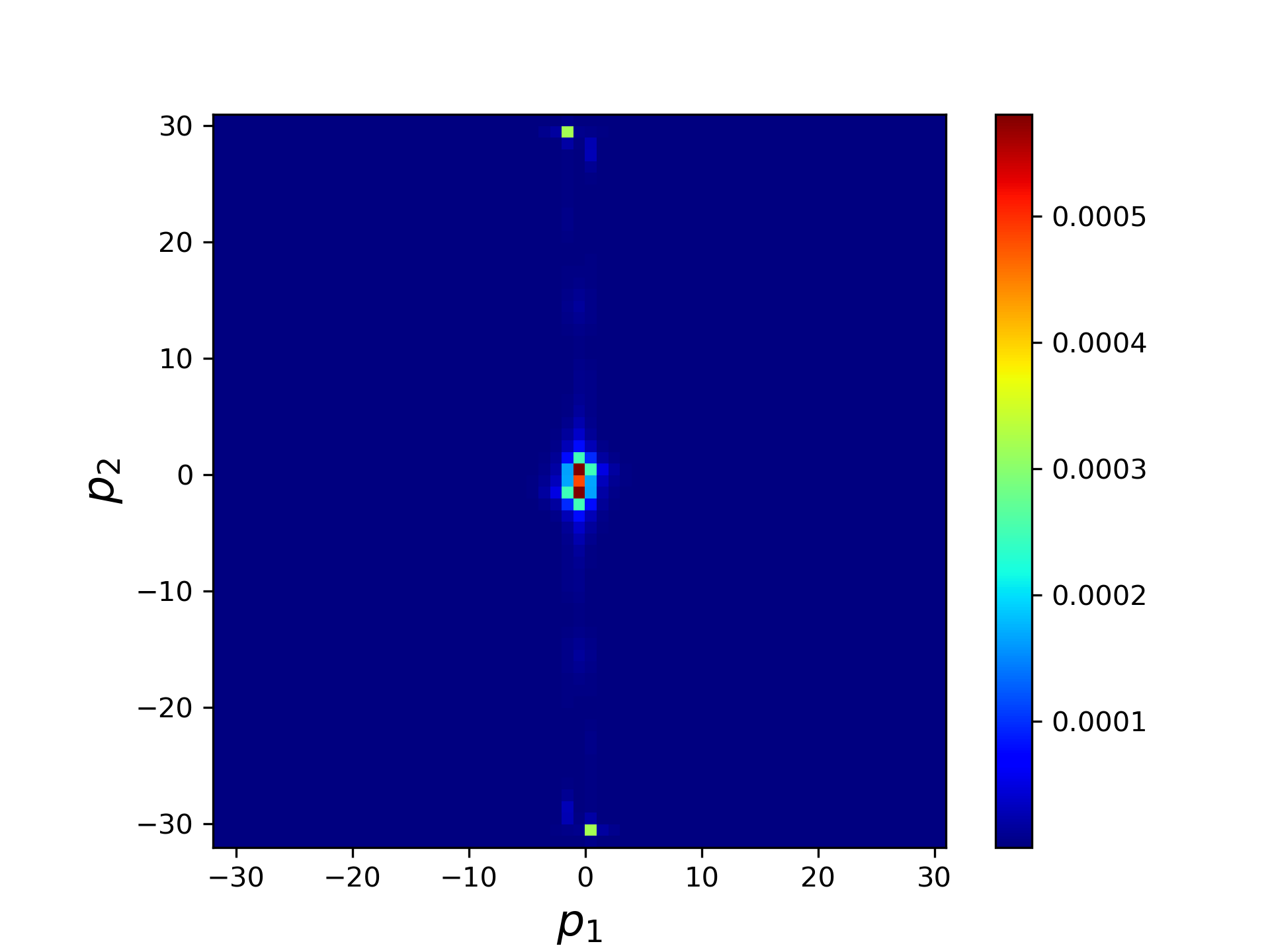}
    }
    \\
    \subfigure[$\xi = 10^{-6}$: convergence factor of $\Phi$]{
        \label{subfig:ACase1_6_a}
        \includegraphics[width=0.27\textwidth]{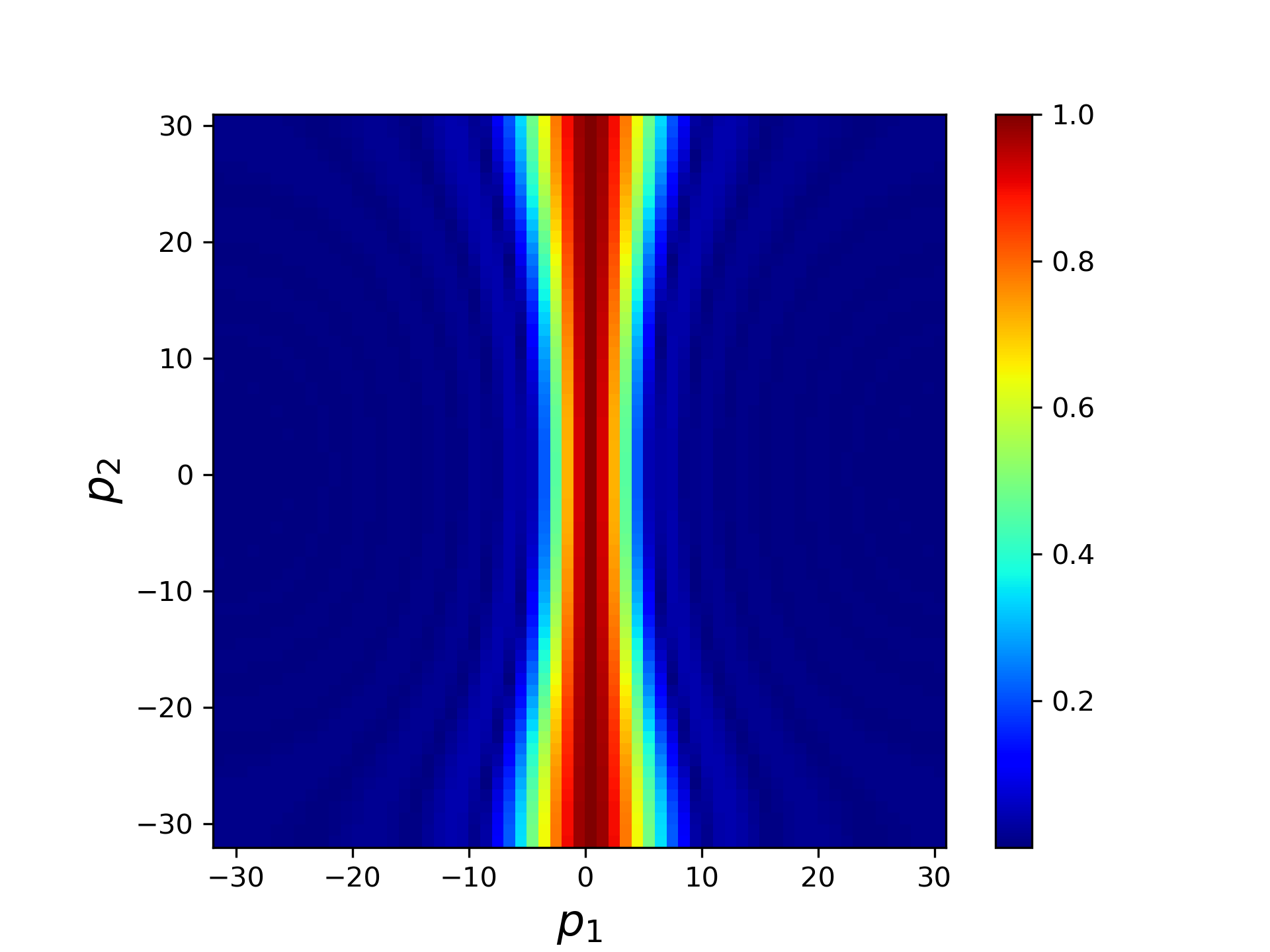}
    }
    \subfigure[$\xi = 10^{-6}$: $\hat{\mathbf{e}}$, before doing corrections]{
            \label{subfig:ACase1_6_b}
            \includegraphics[width=0.27\textwidth]{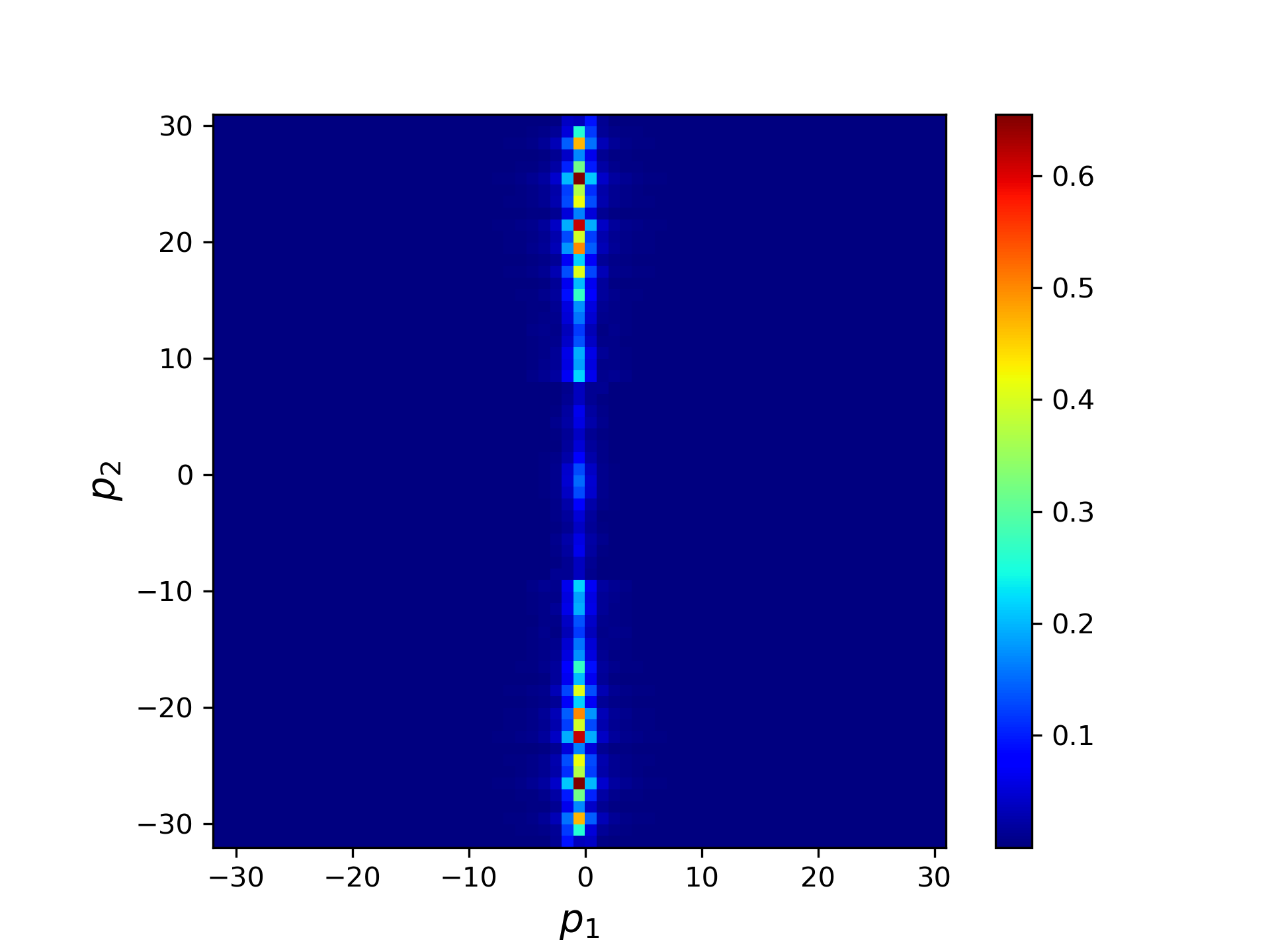}
    }
    \subfigure[$\xi = 10^{-6}$: learned $\hat{\mathbf{e}}$]{
            \label{subfig:ACase1_6_c}
            \includegraphics[width=0.27\textwidth]{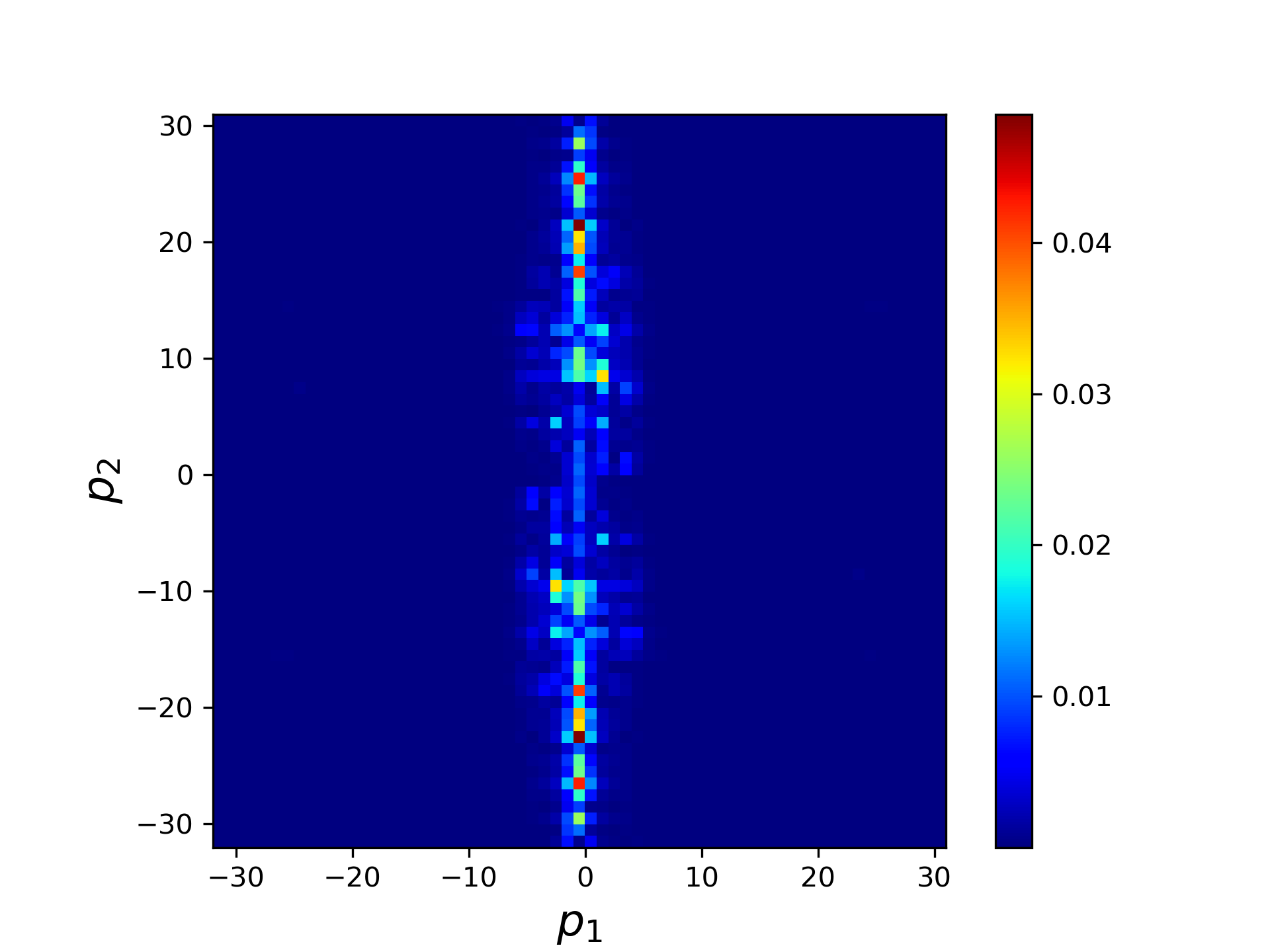}
    }
    \caption{The distribution of convergence factor when $\xi = 10^{-1}, 10^{-6}$. The first column displays the convergence factor of Cheby-semi ($m=10$). The second column shows the error distribution in the frequency space before correction. And the third column shows the error distribution in the frequency space learned by $\mathcal{H}$.}
    \label{fig:Ani_Case1}
\end{figure}

We use $\xi = 10 ^{-1}, 10^{-6}, n=64$ and Cheby-semi ($m=10$) as examples to illustrate the error that $\mathcal{H}$ learned. Figure~\ref{subfig:ACase1_2_a} shows the convergence factor obtained by LFA of Cheby-semi ($m=10$) for solving the system with $\xi = 10 ^{-1}, \theta=0$. It can effectively eliminate the error components except for low-frequency. Figure~\ref{subfig:ACase1_2_b} shows the distribution of error before correction in frequency space. The result is consistent with the guidance of LFA, i.e., the error concentrated in the low-frequency modes. Figure~\ref{subfig:ACase1_2_c} shows the distribution of error learned by $\mathcal{H}$ in the frequency space at this time. Its distribution is largely similar to that of Figure~\ref{subfig:ACase1_2_b}. Figures\,\ref{fig:Ani_Case1}(d-f) show the corresponding situation for $\xi = 10 ^{-6}, \theta=0$. In this case, Cheby-semi ($m=10$) is unable to eliminate the error along the $y$ direction. However, $\mathcal{H}$ is still capable of learning.

\subsubsection{Case 2: Generalization ability of anisotropic direction with fixed strength}
We randomly sample $20$ parameters $\theta$ according to the distribution $\theta \sim U[0,\pi]$ with fixed $\xi = 10^{-6}$. The training and testing data are generated in a similar manner as Section\,\ref{subsec:41} does. Table~\ref{tab:ani_case2} shows test results. It can be seen that whether $\Phi$ is Jacobi or Krylov, FNS can maintain robust performance in all situations, while the line-smoother is not available for these cases.
\begin{table}[h]
  \caption{The mean and standard deviation of the number of iterations required to achieve the stopping criterion over all tests for Eq.\,\eqref{eq:gexiangyixing} case 2.}
  \label{tab:ani_case2}
  \resizebox{\textwidth}{!}{
  \begin{tabular}{@{}lcccccccc@{}}
      \toprule
      $\theta$    & $0.1\pi$         & $0.2\pi$         & $0.3\pi$         & $0.4\pi$         & $0.6\pi$         & $0.7\pi$         & $0.8\pi$         & $0.9\pi$         \\ \midrule
      FNS(Jacobi) & $300.2\pm 23.95$ & $252.4\pm 34.81$ & $269.3\pm 36.70$ & $356.4\pm 39.69$ & $338.4\pm 32.07$ & $265.0\pm 29.96$ & $266.5\pm 25.07$ & $316.7\pm 17.26$ \\ \midrule
      FNS(Krylov) & $58.4\pm 4.45$   & $46.5\pm 4.84$   & $45.1 \pm 2.30 $ & $64.0\pm 7.01$   & $54.4\pm 5.75$   & $41.3\pm 7.11$   & $43.0\pm 2.83 $  & $60.3\pm 3.93$   \\ \bottomrule
      \end{tabular}}
\end{table}

\begin{figure}[!htbp]
    \centering
    \includegraphics[width=13.5cm]{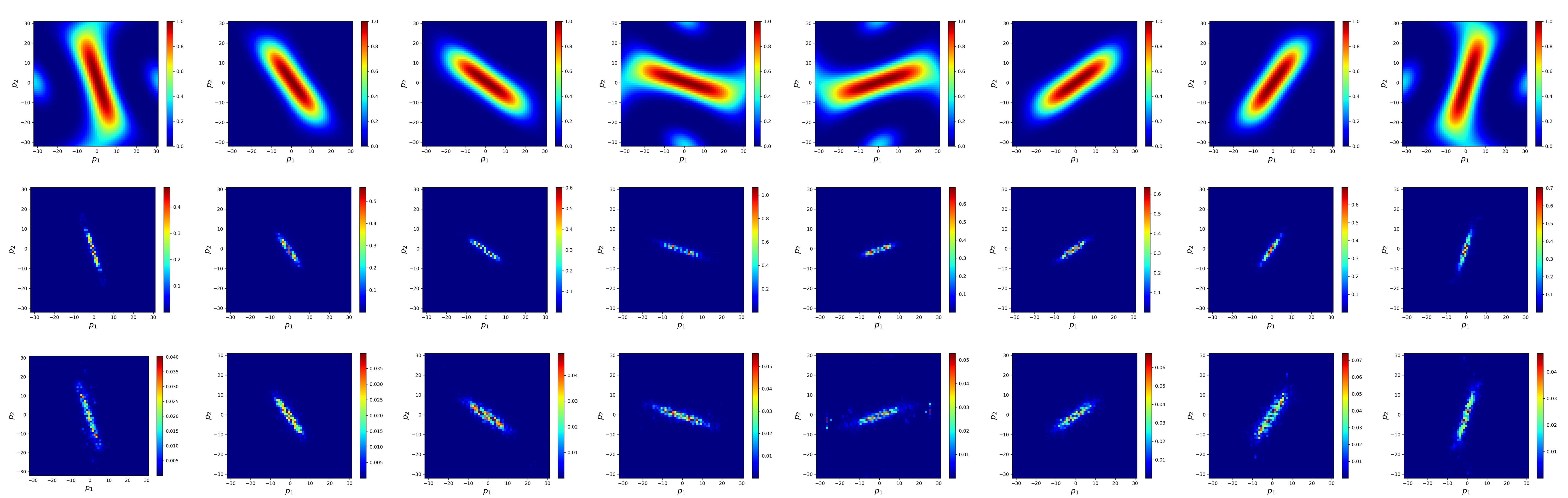}
    \caption{The case 2 of Eq.\,\eqref{eq:gexiangyixing} for different anisotropic direction $\theta$. The first row represents the convergence factor of $\Phi$. The second row represents the error distribution in the frequency space before to correction. The third row represents the learned error by $\mathcal{H}$.}
    \label{fig:Ani_Case2}
\end{figure}

Take $\theta=j\pi/10 , j=1,\ldots,4,6,\ldots, 9$, $\xi=10^{-6}, n=64$ and $\Phi$ is the weighted Jacobi method with weight $2/3$, Figure~\ref{fig:Ani_Case2} shows the test results.
The first row shows the weighted Jacobi method's convergence factor $\mu_\text{loc}$ for each $\theta$ which is computed by LFA.
The region of $\mu_\text{loc}\sim 1$ means that the error is difficult to eliminate. It can be found that these error components are distributed along the anisotropic direction. The second row depicts the error distribution in frequency space before correction, which is consistent with the results obtained by the LFA. The third row shows the distribution of the error learned by $\mathcal{H}$. It can be seen that $\mathcal{H}$ can automatically learn the error components that $\Phi$ is difficult to eliminate. Line-smoother is not feasible for these $\theta$.

\subsection{Convection-diffusion equation}
Consider the convection-diffusion equation
\begin{equation}
    \begin{cases}-\varepsilon \Delta u+u_{x}+u_{y}=0, & \Omega=(0,1)^{2}, \\ u=0, & x=0,0 \leq y<1 \text { and } y=0,0 \leq x<1, \\ u=1, &x=1,0 \leq y<1 \text { and } y=1,0 \leq x \leq 1,\end{cases}
    \label{eq:conv-diff}
\end{equation}
We use central difference method to discretize \eqref{eq:conv-diff} on a uniform mesh with spatial size $h$ in both $x$ and $y$ directions which yields a non-symmetric stencil
\begin{equation}
    \frac{1}{h^{2}}\left[\begin{array}{ccc}
    & h / 2-\varepsilon & \\
    -h / 2-\varepsilon & 4 \varepsilon & h / 2-\varepsilon \\
    & -h / 2-\varepsilon &
    \end{array}\right].
\end{equation}

For the requirement of stability, the central difference scheme needs to satisfy the Peclet condition
\begin{equation}
    P e:=\frac{h}{\varepsilon} \max (|a|,|b|) \leq 2,
\end{equation}
which means that the central difference method cannot approximate the PDE solution when $\varepsilon$ is extremely small. However, here we only take into account the solver of the linear system. Thus we continue to use this discretization method to demonstrate the performance of FNS. 
In the next experiments, we will explore the diffusion- and convection-dominant cases, respectively.

\subsubsection{Case 1: $\varepsilon \in (0.01, 1)$}
\begin{figure}[!htbp]
    \centering
    \subfigure[Convergence factor of one times weighted Jacobi ($\omega=4/5$)]{\label{subfig:JLFA}\includegraphics[width=0.23\textwidth]{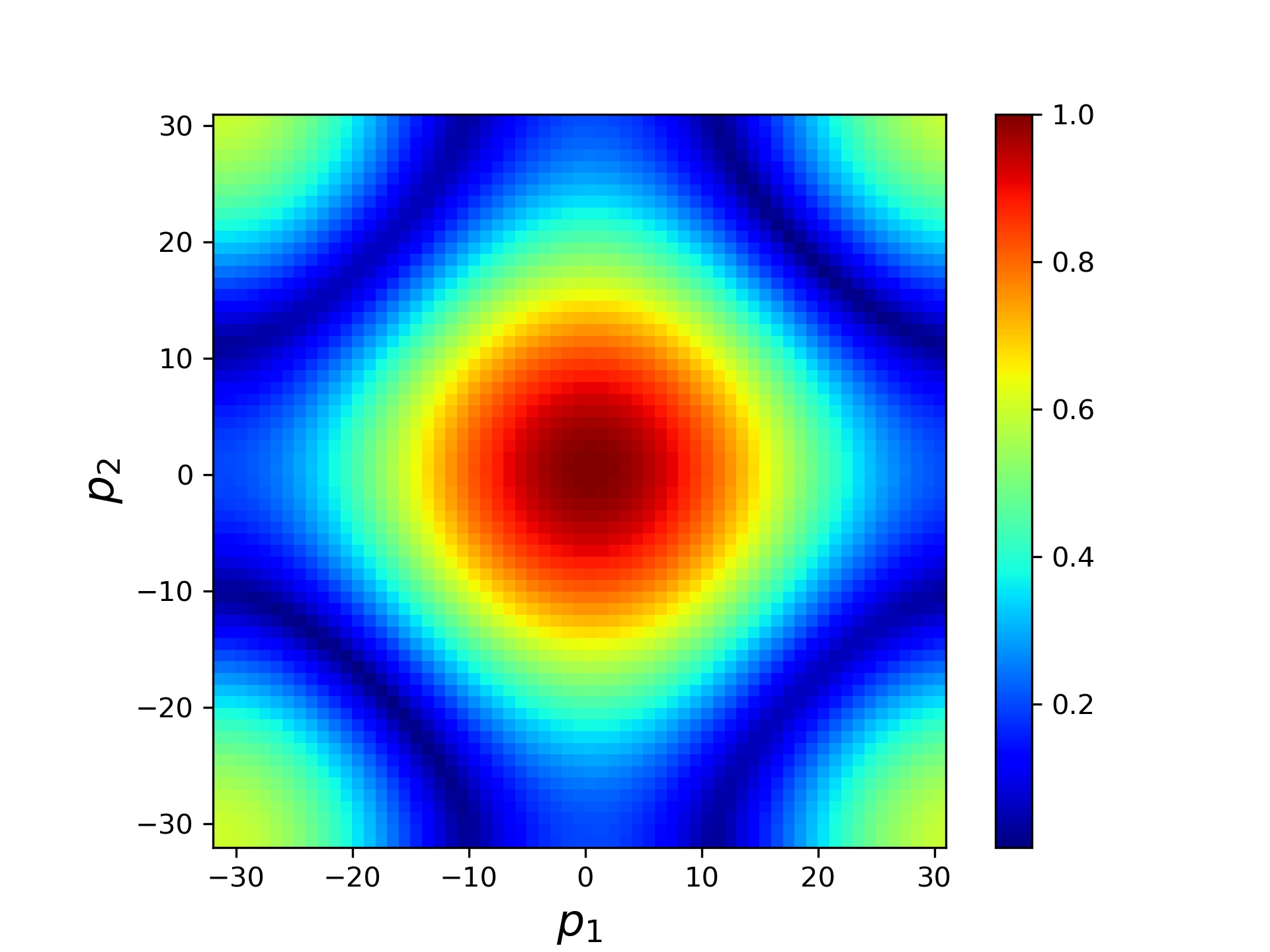}}
    \subfigure[Two times]{\label{subfig:b}\includegraphics[width=0.23\textwidth]{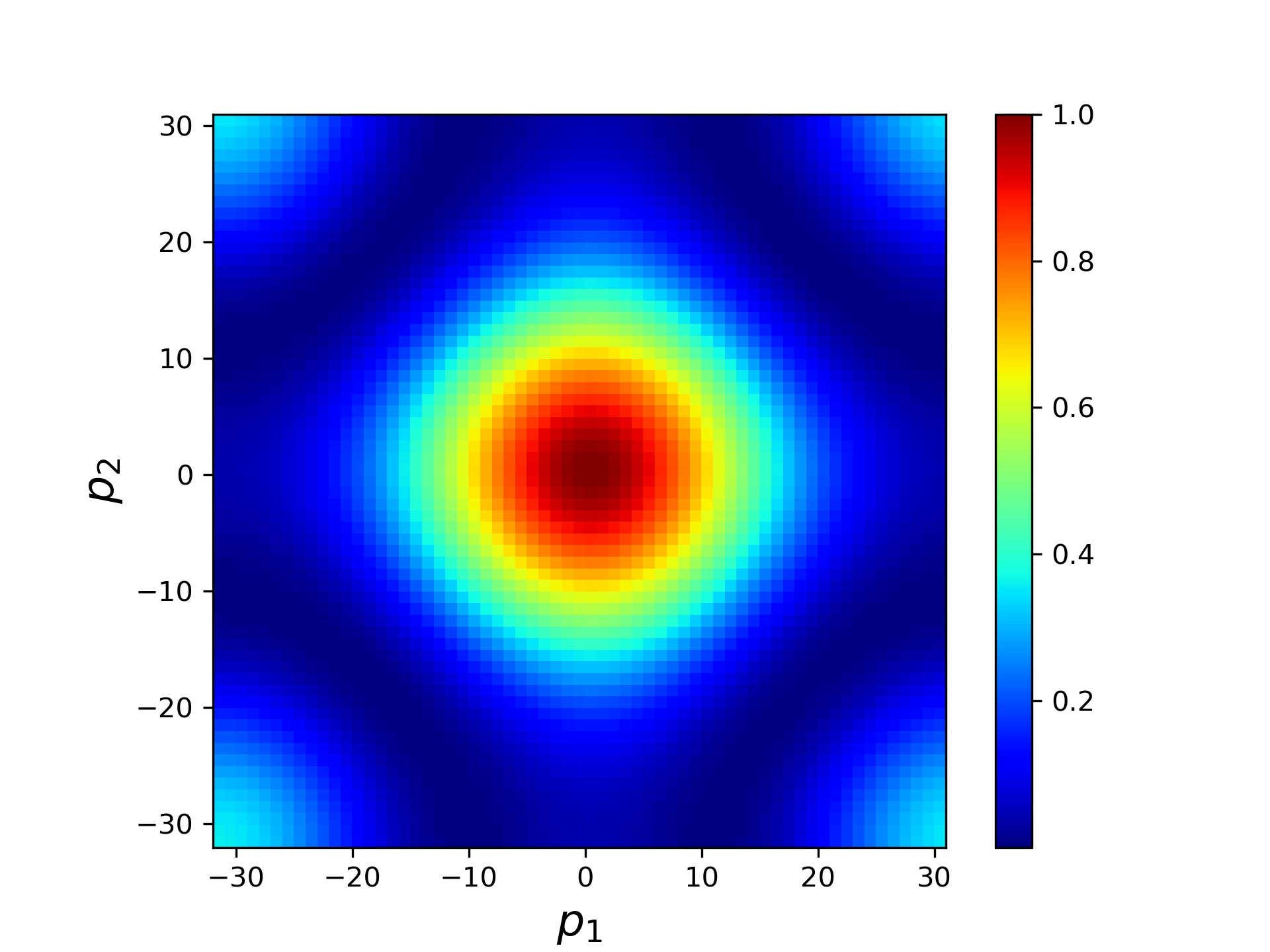}}
    \subfigure[Three times]{\includegraphics[width=0.23\textwidth]{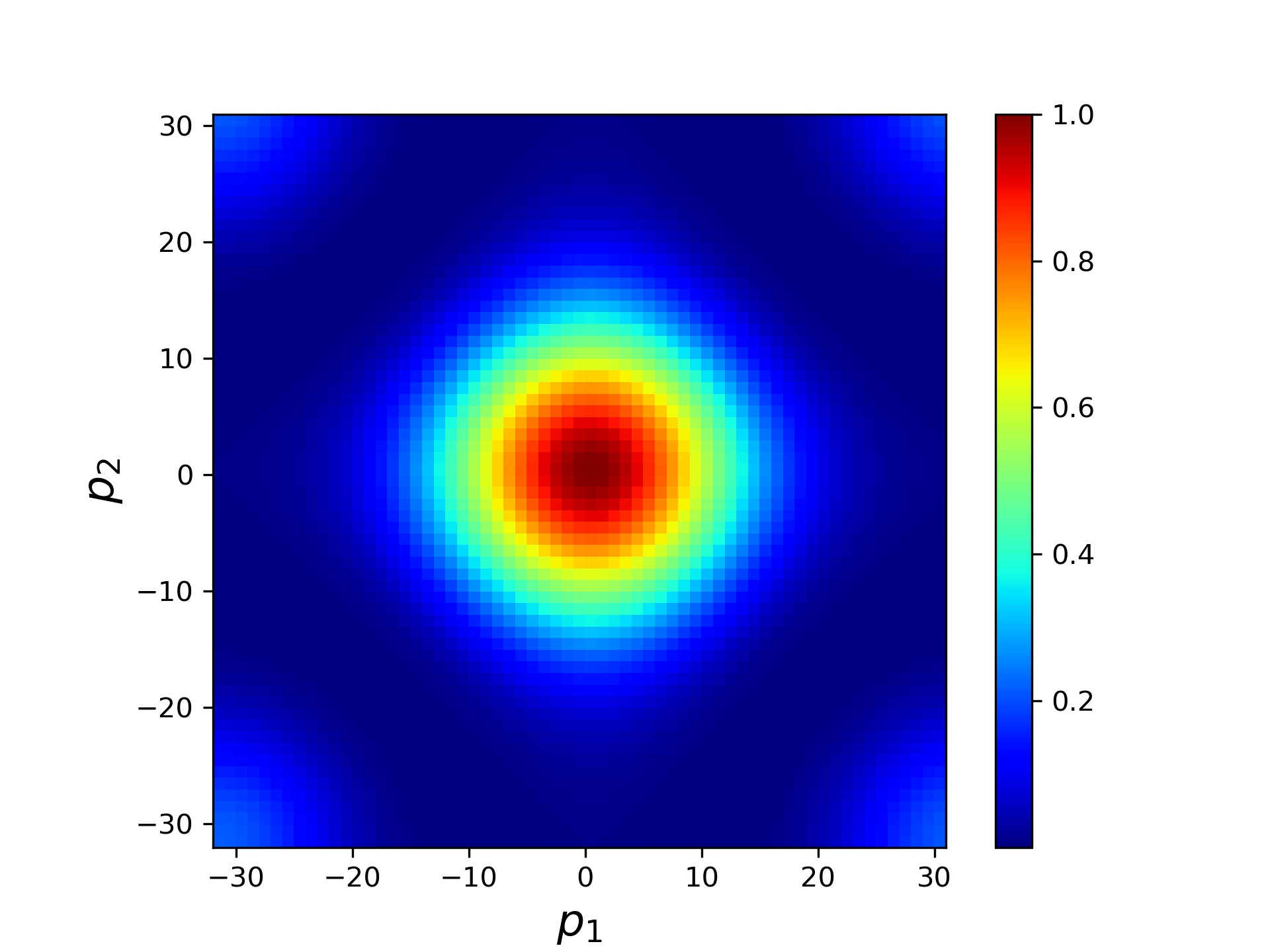}}\\
    \subfigure[Four times]{\includegraphics[width=0.23\textwidth]{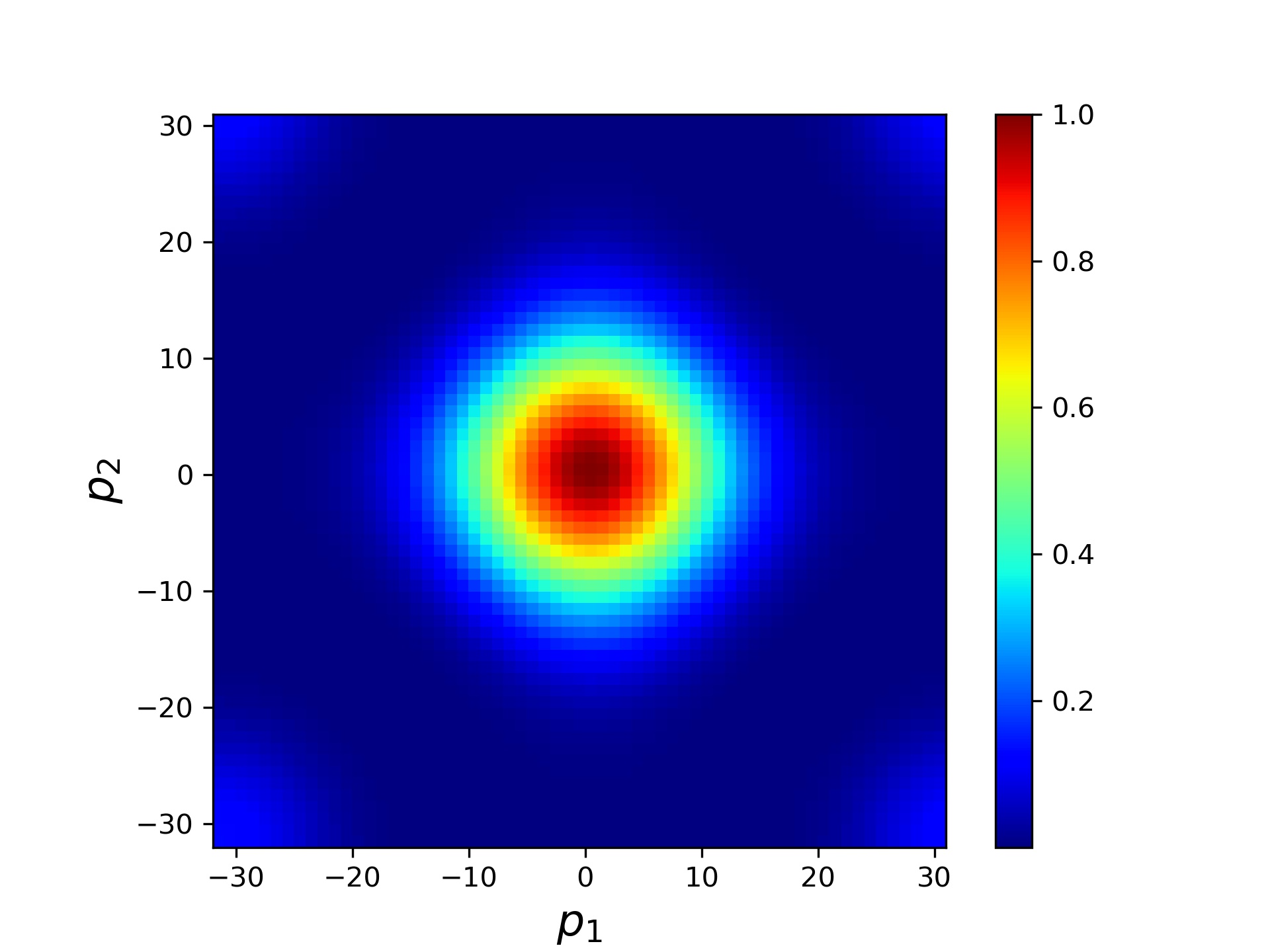}}
    \subfigure[Five times]{\label{subfig:e}\includegraphics[width=0.23\textwidth]{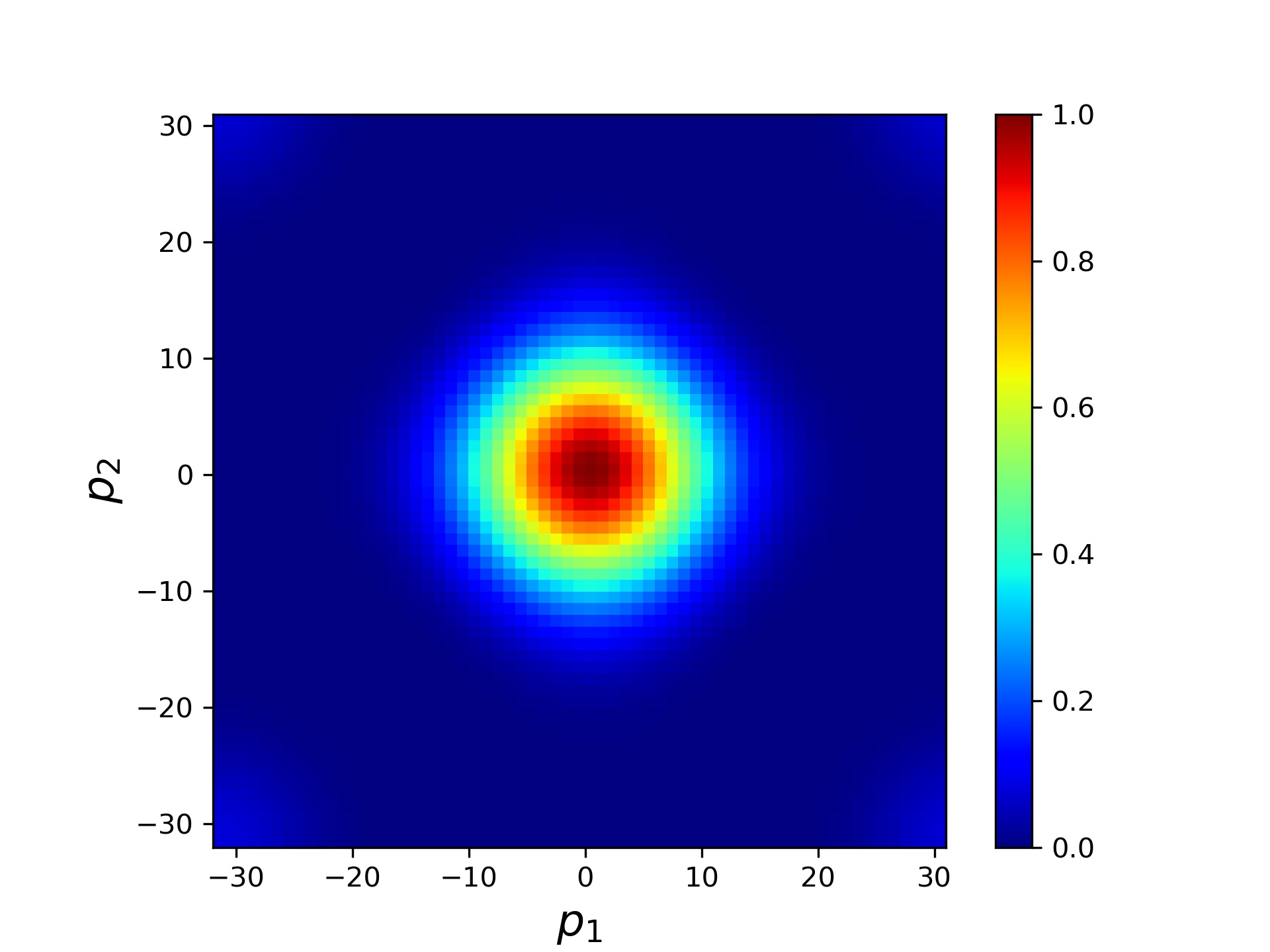}}
    \subfigure[Learned $\hat{\mathbf{e}}$]{\label{subfig:Lehat}\includegraphics[width=0.23\textwidth]{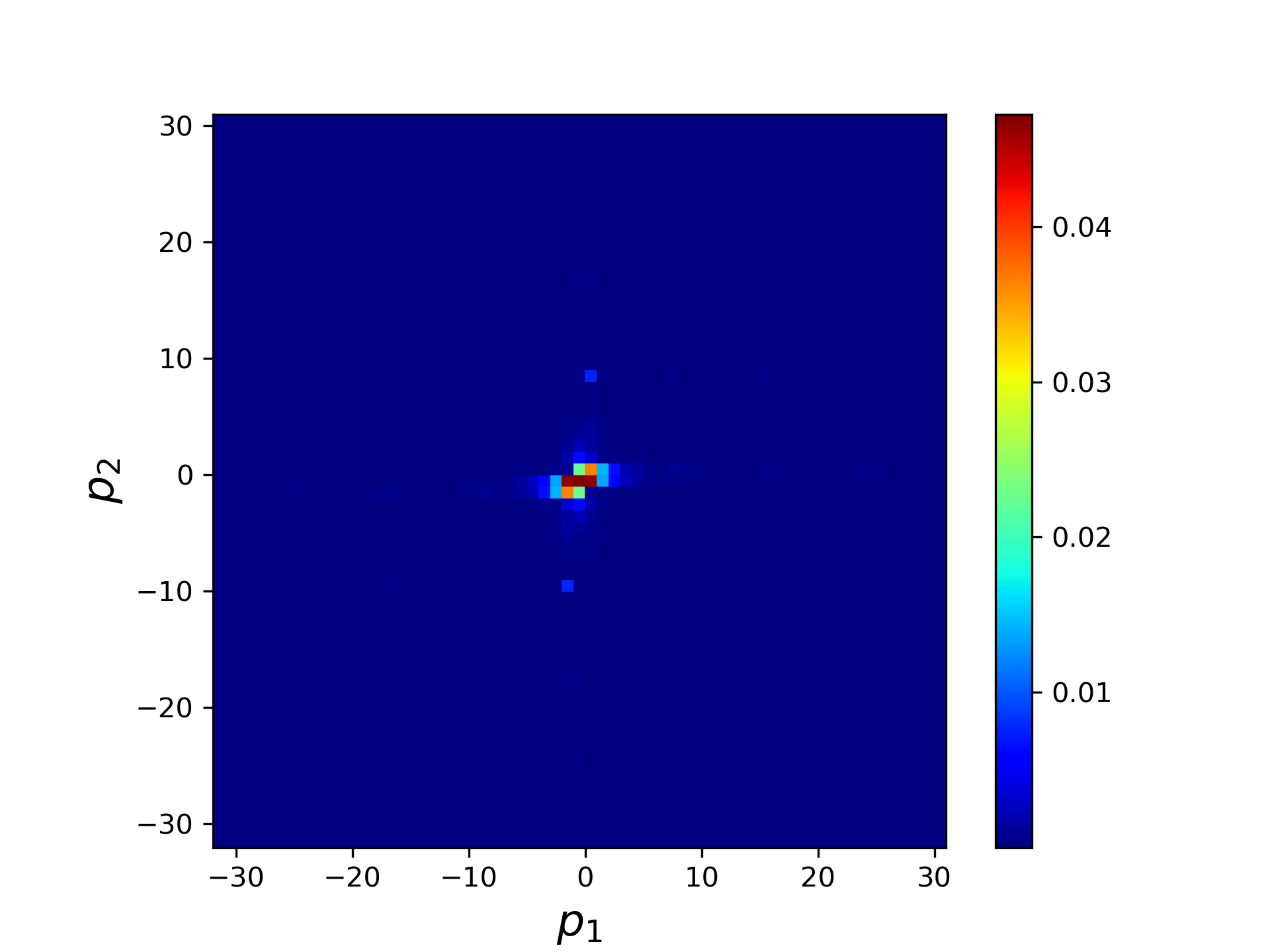}}
    \caption{The distribution of convergence factor of 
    five-times consecutive weighted Jacobi method. The last plot shows the distribution of errors learned by $\mathcal{H}$ in frequency space.}
\end{figure}

We utilize weighted Jacobi method as $\Phi$ in this case. Taking $\varepsilon=0.1, h = 1/64$ as an example, Figure~\ref {subfig:JLFA} illustrates the convergence factor obtained by LFA of the weighted Jacobi method ($\omega=4/5$) for solving system. We use five-times consecutive weighted Jacobi method as $\Phi$. Figure~\ref{subfig:b}-\ref{subfig:e} show that such $\Phi$ is a good smoother. Figure~\ref{subfig:Lehat} shows the distribution of  error learned by $\mathcal{H}$ in the frequency space. It can be observed that this is essentially complementary to $\Phi$.

\begin{figure}[!htbp]
    \centering
    \subfigure[$\varepsilon=0.5$]{\includegraphics[width=0.3\textwidth]{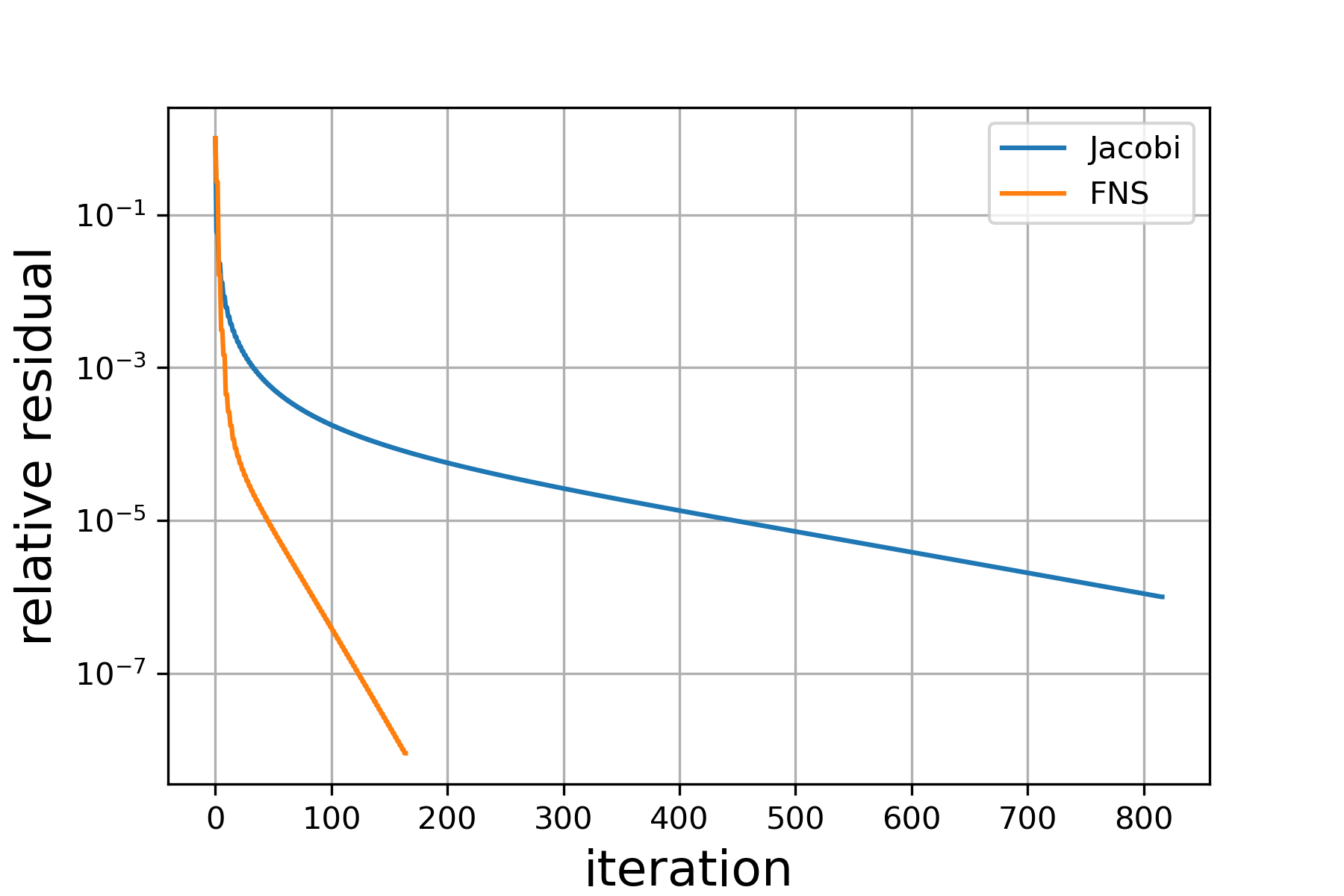}}
    \subfigure[$\varepsilon=0.1$]{\includegraphics[width=0.3\textwidth]{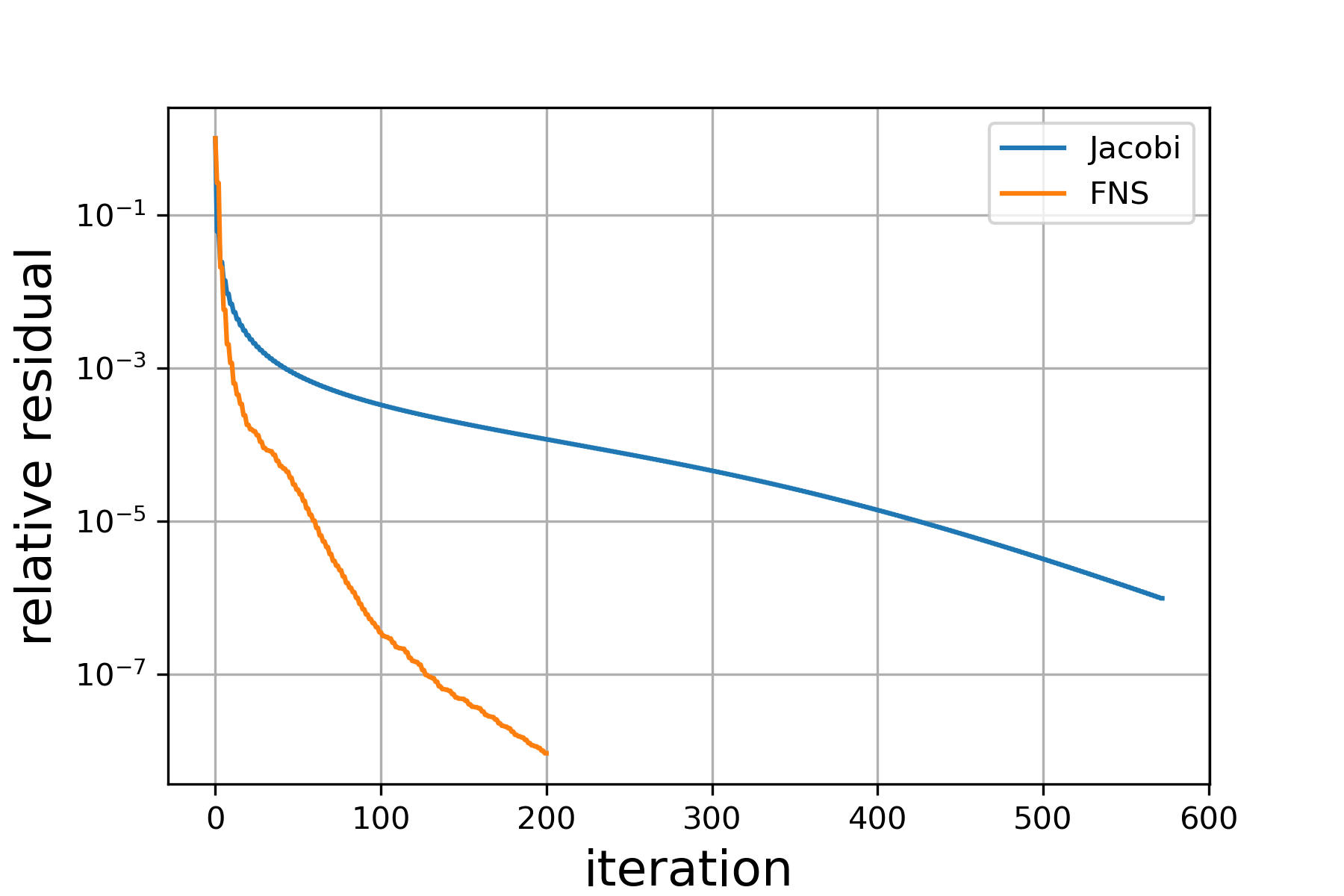}}
    \subfigure[$\varepsilon=0.05$]{\includegraphics[width=0.3\textwidth]{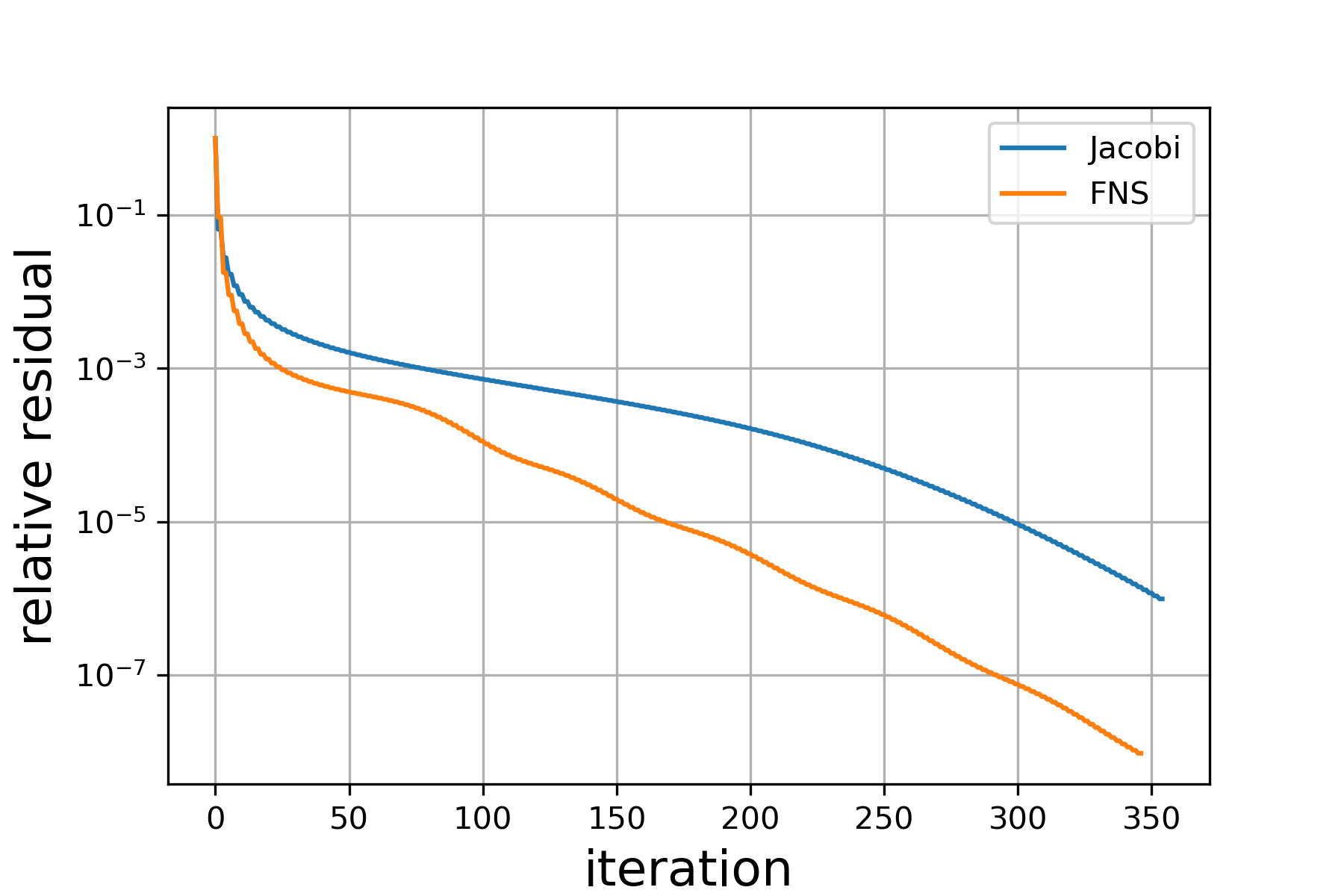}}
    \caption{The relative residual with the FNS and weighted Jacobi method, where Jacobi denotes five-times consecutive weighted Jacobi method.}
    \label{fig:residual_ad}
\end{figure}

Figure~\ref{fig:residual_ad} uses $\varepsilon=0.5, 0.1, 0.05$ as examples to show the relative residual of FNS and weighted Jacobi method. It can be seen that  FNS has acceleration and the weighted Jacobi method ramps up as $\varepsilon$ decreases. This is due to the fact that when $\varepsilon$ declines, the diagonal element $4\varepsilon$ becomes small, and the weight along the gradient direction increases. However, Jacobi and other gradient descent algorithms will diverge as long as $\varepsilon$ continues to decrease unless the weight is drastically lowered.

\subsubsection{Case 2: $\varepsilon \in [10^{-6}, 10^{-3}]$}
\begin{figure}[!htbp]
    \centering
    \subfigure[$\varepsilon=10^{-2}$]{\includegraphics[width=0.3\textwidth]{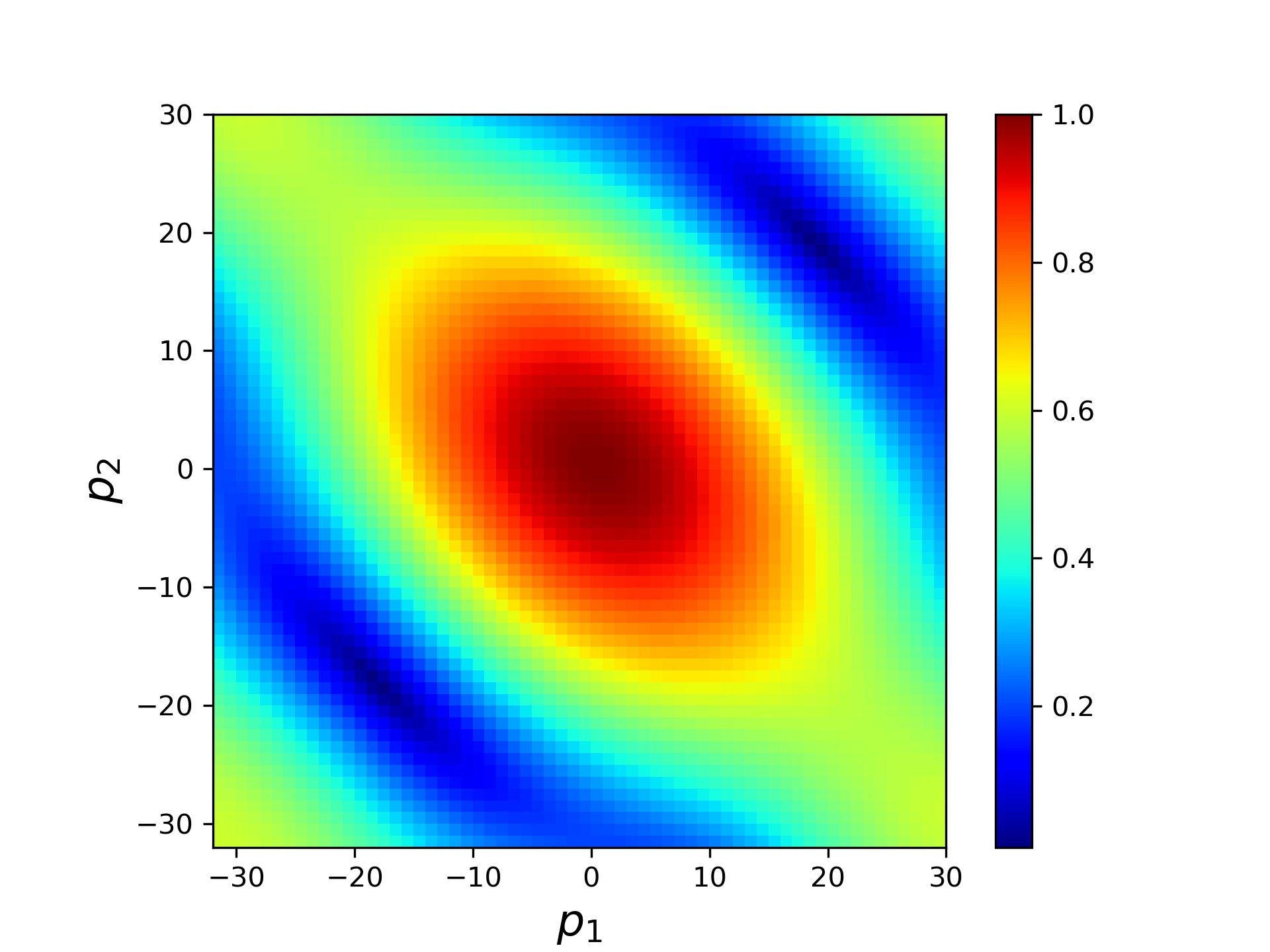}}
    \subfigure[$\varepsilon=10^{-3}$]{\includegraphics[width=0.3\textwidth]{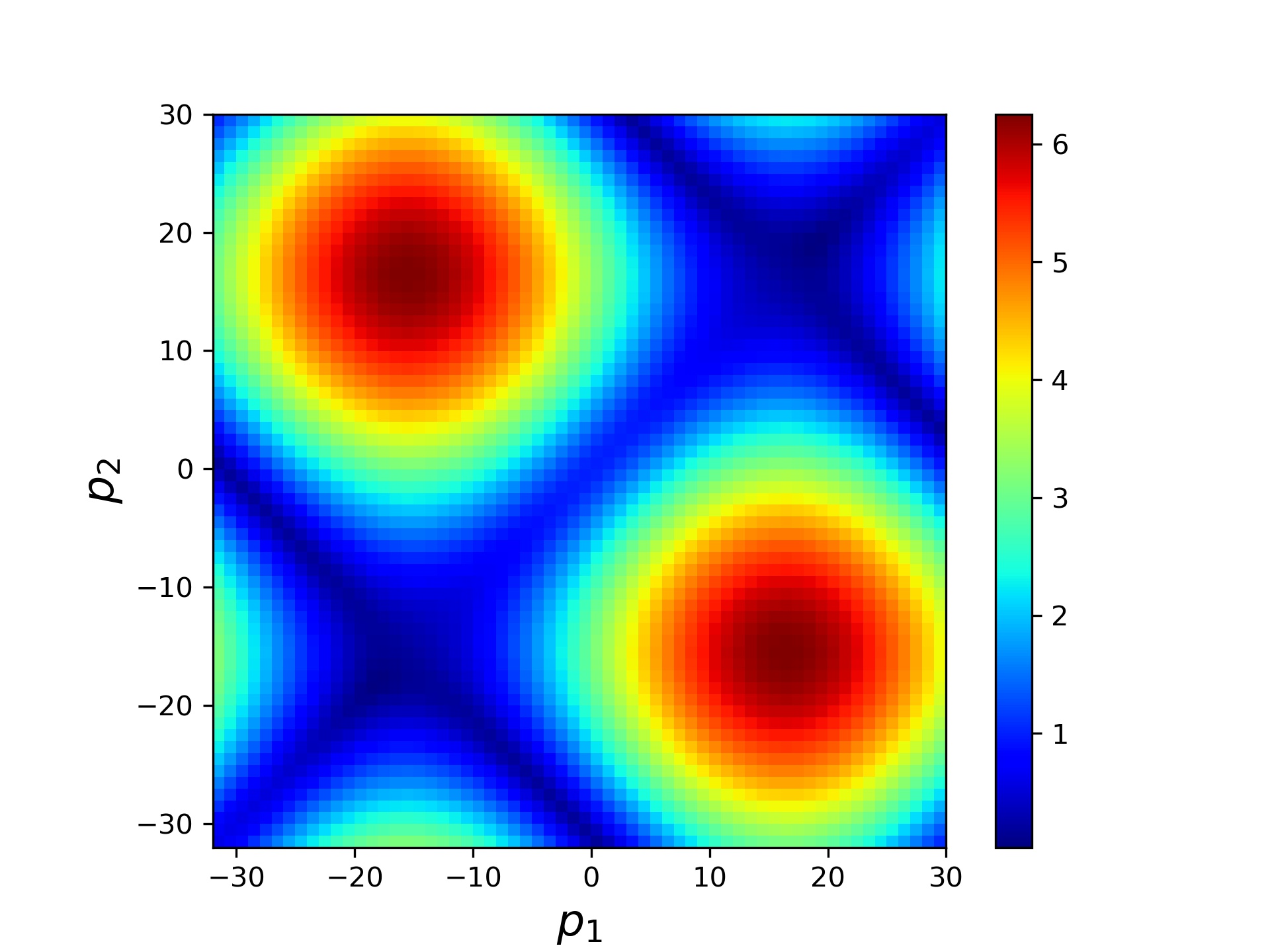}}
    \subfigure[$\varepsilon=10^{-6}$]{\includegraphics[width=0.3\textwidth]{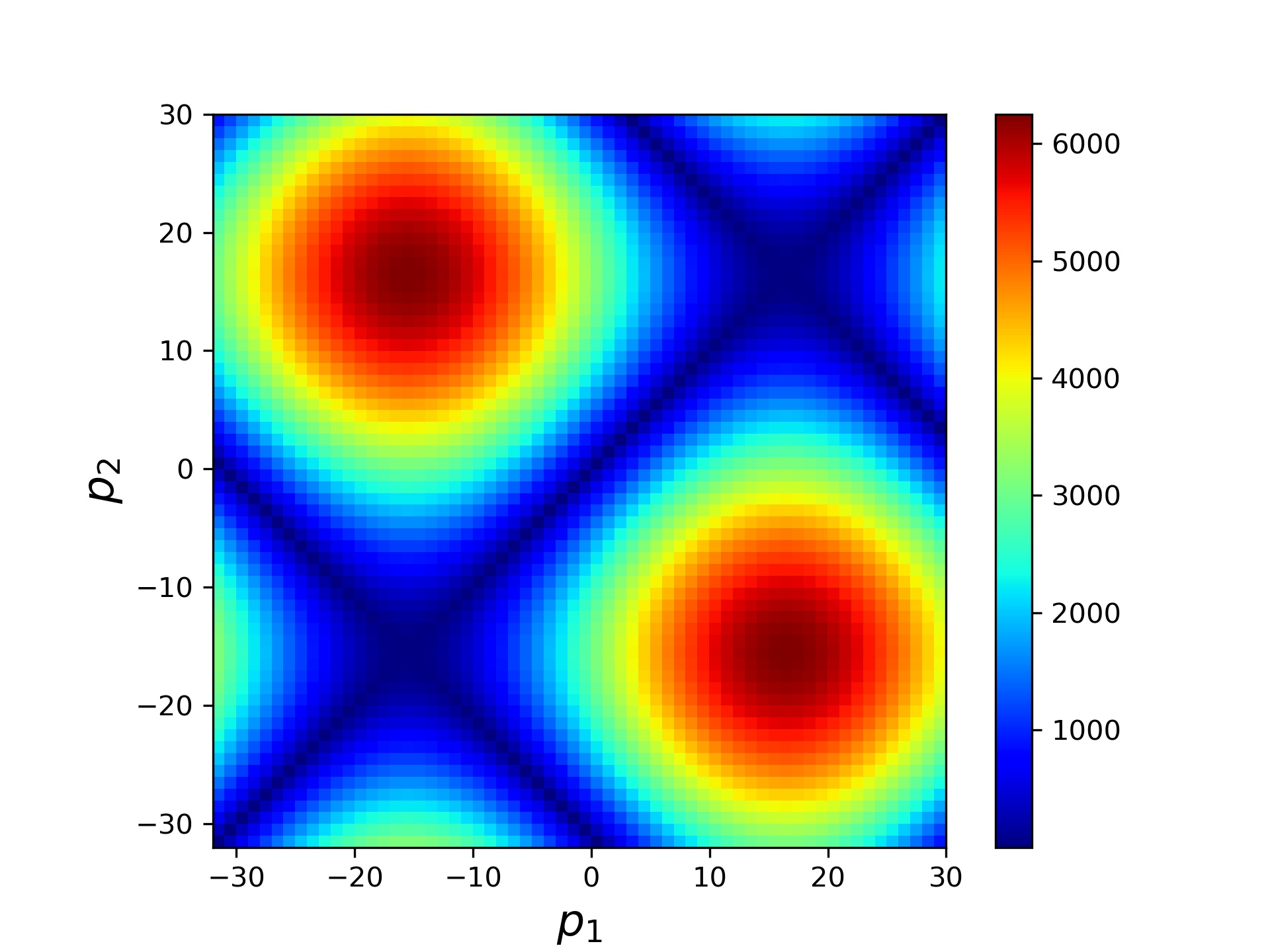}}
    \caption{The distribution of convergence factor for weighted Jacobi method ($\omega=4/5$) when solving the system corresponding to $\varepsilon=10^{-2},10^{-3},10^{-6}$.}
    \label{fig:LFA_ad_Jacobi}
\end{figure}

In this case, since the diagonal elements of the discrete system are notably less than the off-diagonal elements, the system is non-symmetric. Many methods such as Jaocbi, CG, and MG (Jacobi) methods might diverge. Figure~\ref{fig:LFA_ad_Jacobi} shows that convergence factor of weighted Jacobi ($\omega=4/5$) for solving the system when $\varepsilon=10^{-2},10^{-3},10 ^{-6}$. It can be seen that when $\varepsilon$ is small, the convergence factor of weighted Jacobi method for most frequency modes is bigger than $1$, which causes this iterative method to diverge and be unsuitable as a smoother.

Consequently, we learn the $\Phi$ in Eq.\,\eqref{eq:simple_iter}, where $\mathbf{B}$ is a two-layer linear CNN with channels 1$\rightarrow$8$\rightarrow$1, and the kernel size is $3\times 3$. The $\Phi$ is trained together with $\mathcal{H}$, and the training hyperparameters are listed in Appendix\,\ref{appen:ad}.
Figure~\ref{subfig:small_diff} illustrates how this learned FNS is able to solve the linear system when $\varepsilon=10^{-3}, 10^{-4}, 10^{-5},10^{-6}$. It is visible that FNS converges soon. 
Figure~\ref{subfig:gmres} shows the change of relative residual for FNS, GMRES, and BiCGSTAB($\ell$) ($\ell=15$) with $\varepsilon=10^{-6}$. It is clear that FNS has the fastest convergence rate.
\begin{figure}[!htbp]
    \centering
    \subfigure[different $\varepsilon$]{\label{subfig:small_diff}\includegraphics[width=0.4\textwidth]{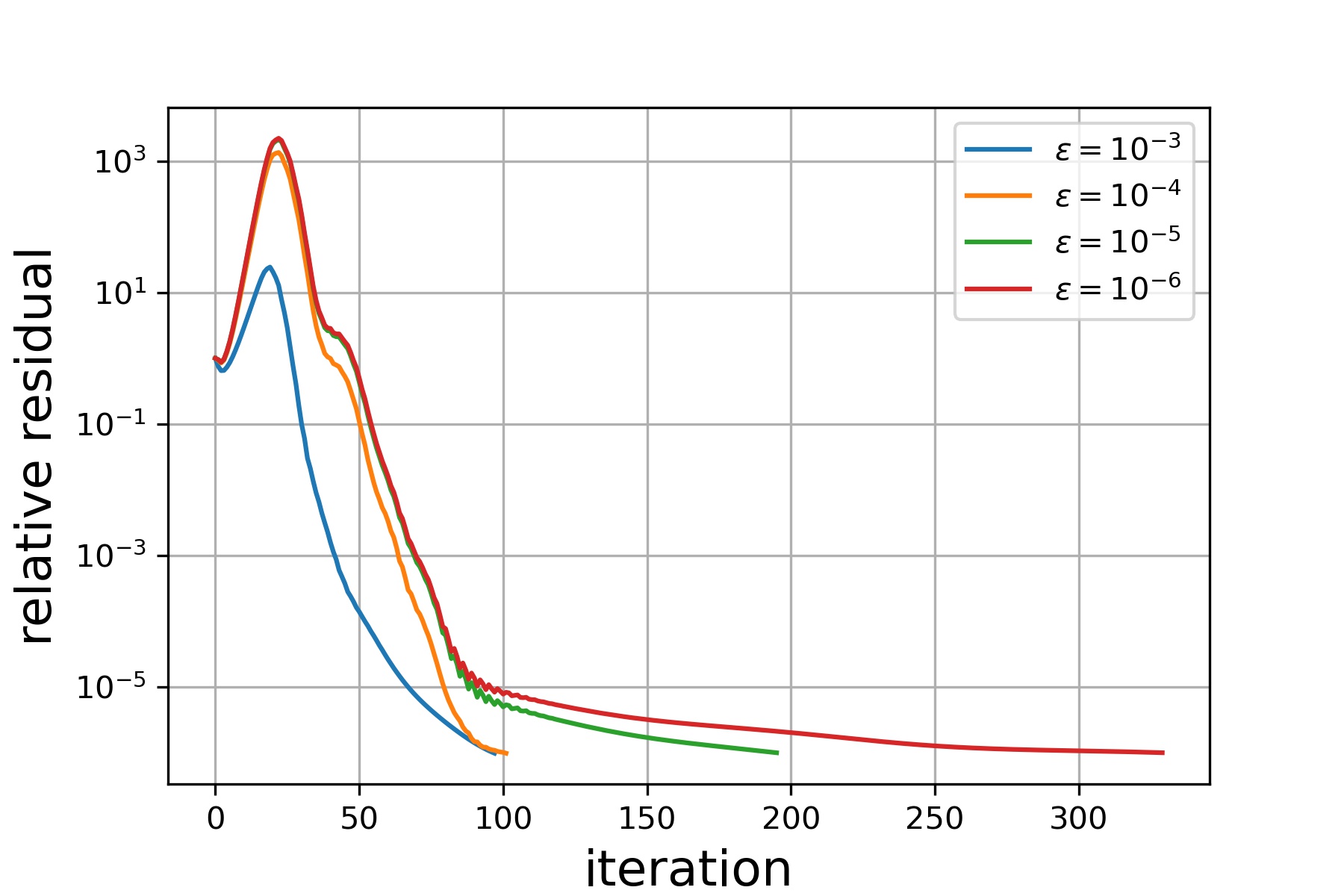}}
    \subfigure[$\varepsilon=10^{-6}$]{\label{subfig:gmres}\includegraphics[width=0.4\textwidth]{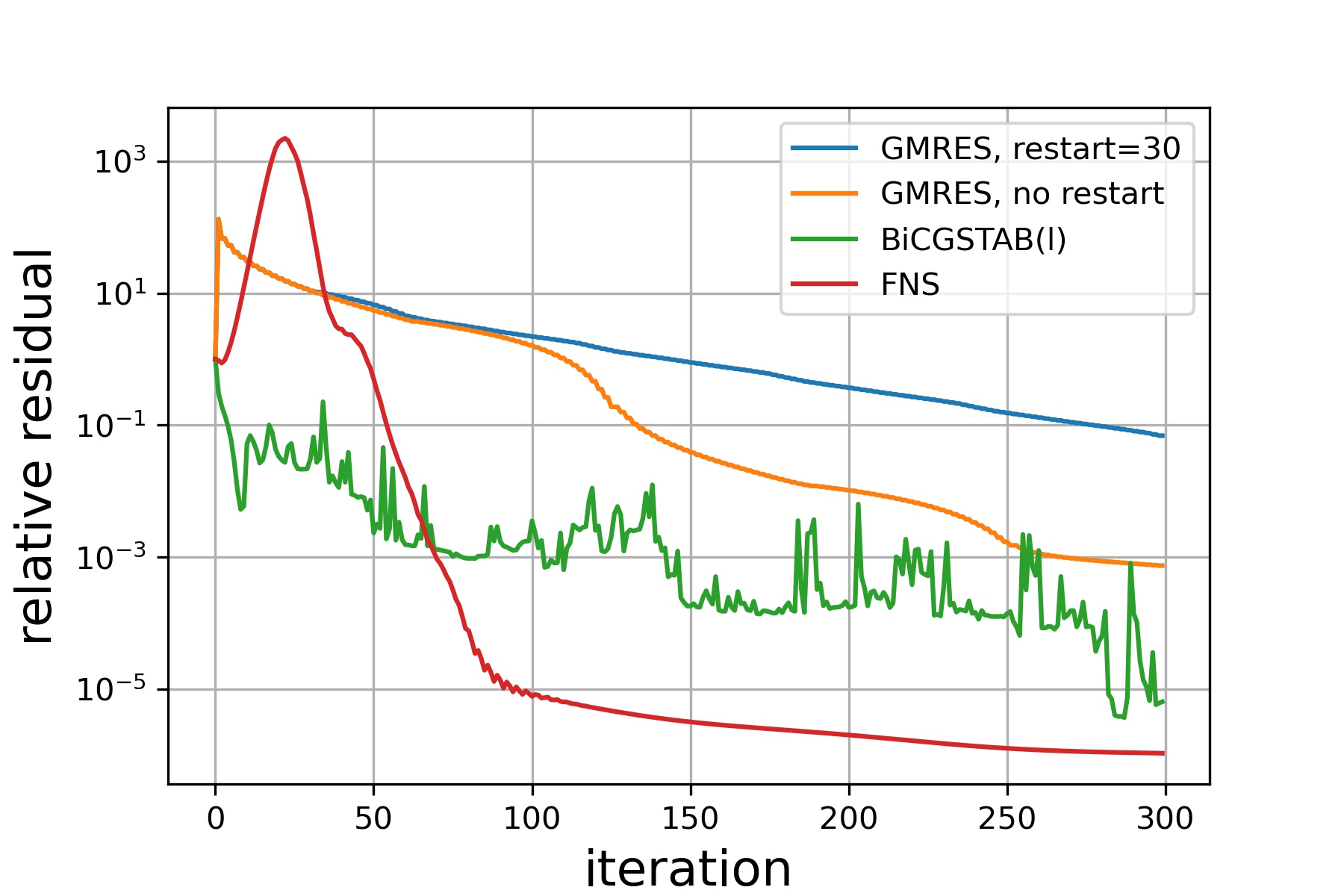}}
    \caption{(a): The change of relative residual with FNS iteration steps for different $\varepsilon$. (b): Comparison of FNS and GMRES, BiCGSTAB($\ell$) when $\varepsilon=10^{-6}$.}
\end{figure}
 
\subsection{Helmholtz equation}
The Helmholtz equation we consider here is 
\begin{equation}
    -\Delta u(x)- \kappa^2 u(x)=g(x), \quad x \in \Omega,
    \label{eq:helm}
\end{equation}
where $\Omega = (0,1)^2$, $\kappa$ is the wavenumber. 
We currently only take into account the zero Dirichlet boundary condition. We use the second order FDM to discretize \eqref{eq:helm} on a uniform mesh with spatial size $h$. The corresponding stencil reads
\begin{equation}
    \frac{1}{h^2}\left[\begin{array}{ccc}
    0 & -1 & 0 \\
    -1 & 4-\kappa^2 h^2 & -1 \\
    0 & -1 & 0
    \end{array}\right].
\end{equation}

We examine the FNS performance at a low wavenumber ($\kappa=25$) and medium wavenumber ($\kappa=125$). For $\kappa=25$, we take $h = 1/64$, and $h = 1/256$ for $\kappa=125$. Take Krylov in \cite{chen2022meta} as $\Phi$, and $g(x)=1$, the training hyperparameters are listed in Appendix\,\ref{appen:helm}.
Figure~\ref{fig:helmholtz} depicts how the relative residual decreased with different solvers.
For $\kappa=25$, FNS performs best for the first $300$ steps, but BiCGSTAB performs better at the end. For $\kappa=125$, FNS outperforms BiCGSTAB. And the GMRES results were too subpar to display in this case.
\begin{figure}[!htbp]
    \centering
    \subfigure[$\kappa=25, h = \frac{1}{64}$]{\label{subfig:small_k}\includegraphics[width=0.4\textwidth]{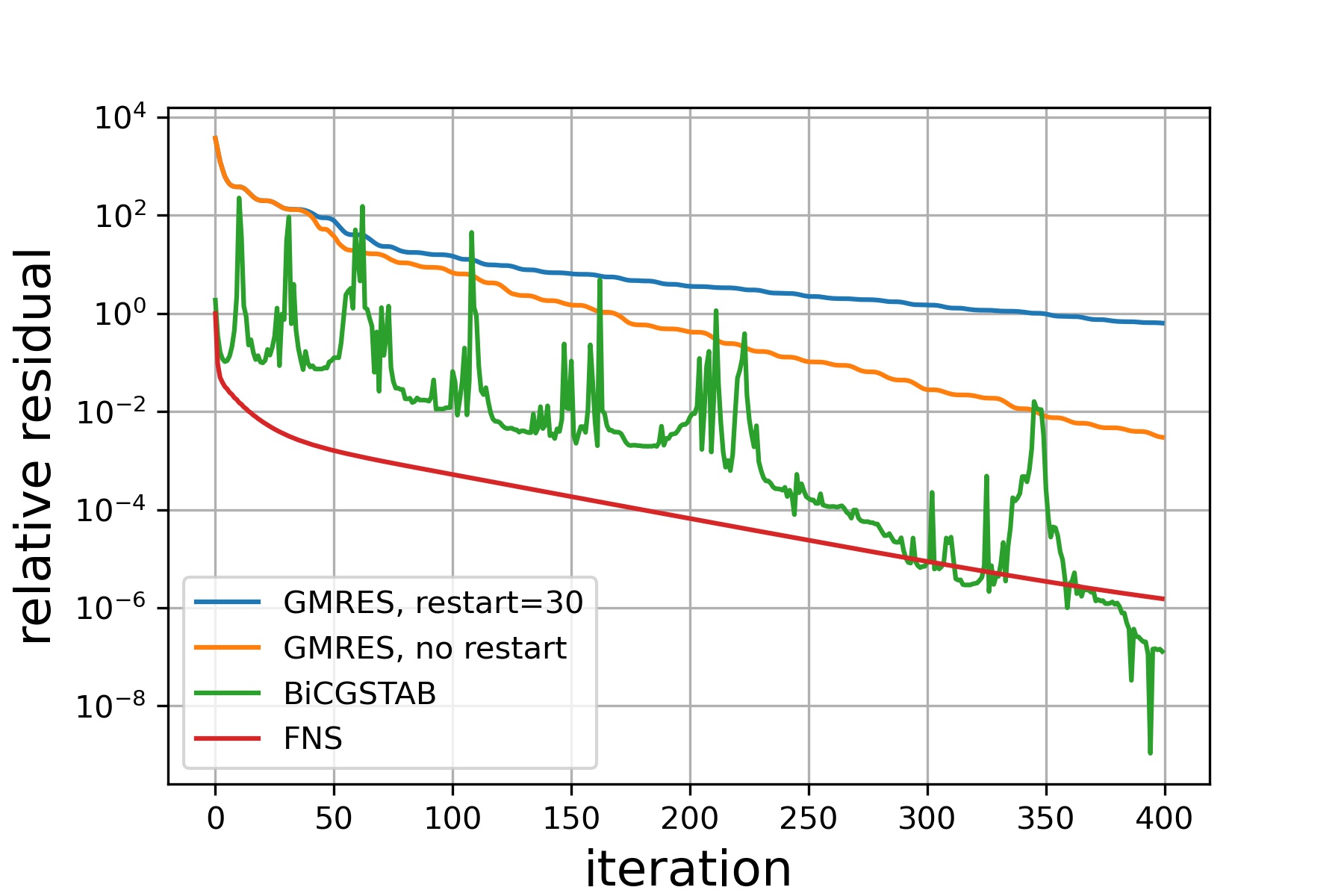}}
    \subfigure[$\kappa=125,h = \frac{1}{256}$]{\label{subfig:mediumk}\includegraphics[width=0.4\textwidth]{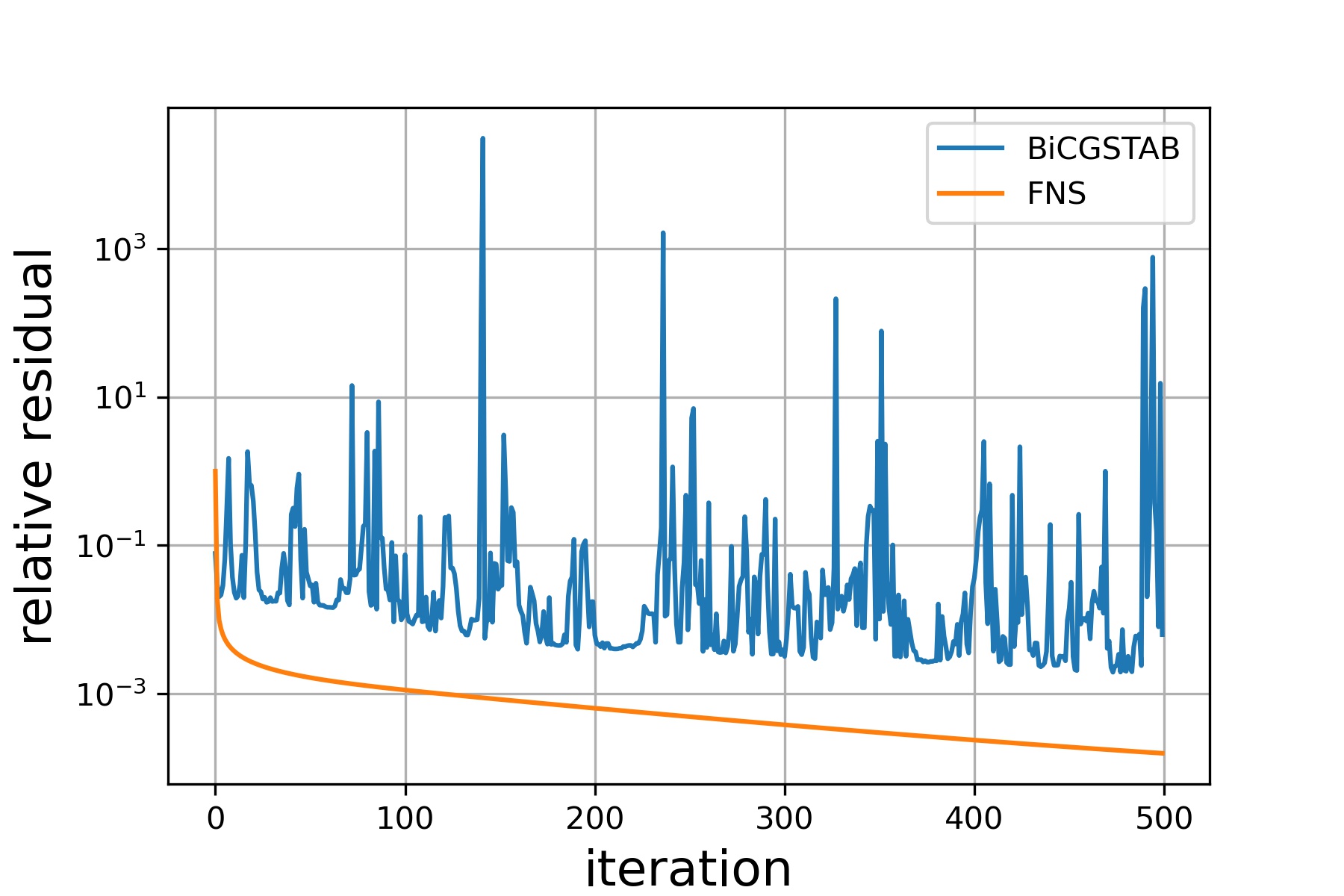}}
    \caption{The change of relative residual for FNS and other solvers when $\kappa=25$ and $\kappa=125$,
    respectively. }
    \label{fig:helmholtz}
\end{figure}

\section{Conclusions and future work}\label{sec:05}
This paper proposes an interpretable FNS to solve large sparse linear systems. It is composed of the stationary method and frequency correction, which are used to eliminate errors in different Fourier modes.
Numerical experiments show that FNS is more effective and robust than other solvers in solving anisotropy equation,  convection-diffusion equation and Helmholtz equation. The core concepts discussed here can be available to a broad range of systems.

There is still a great deal of work to do. Firstly, we only considered uniform mesh in this paper. We will generalize FNS to non-uniform grids, by exploiting tools from geometric deep learning, such as graph neural network and graph Fourier transform. 
Secondly, as we discussed previously, the stationary method converges slowly or diverges in some situations, which prompts researchers to approximate solutions in other transform space. This is true for almost advanced iterative methods, including MG, Krylov subspace methods and FNS. This specified space, however, may not always be the best choice. In the future, we will investigate additional potential transforms, such as Chebyshev, Legendre transforms, and even learnable transforms based on data.

\section{Training hyperparameters}
\subsection{Anisotropy equation}\label{appen:ani}
\begin{table}[h]
    \centering
    \caption{Training hyperparameters for anisotropic equation.}
    \label{tab:ani_hyper_case1}
    \begin{tabular}{@{}lccccc@{}}
        \toprule
        & learning rate & batch size & $K$  & xavier init & grad clip \\ \midrule
FNS(Cheby-semi) & $10^{-4}$     & $100$      & $10$ & $10^{-2}$   & false     \\ 
FNS(Jacobi)     & $10^{-4}$     & $100$      & $10$ & $10^{-2}$   & false     \\
FNS(Krylov)     & $10^{-4}$     & $100$      & $10$ & $10^{-2}$   & false    \\ \bottomrule
\end{tabular}
\end{table}
\begin{table}[h]
    \centering
    \caption{HyperNN architecture parameters for Anisotropic equation.  Notations in ConvTranspose2d are: i: in channels; o: out channels; k: kernel size; s: stride; p: padding.}
    \label{tab:ani_hyper_nn}
    \begin{tabular}{@{}l@{}}
        \toprule
        ConvTranspose2d(i=1,o=4,k=3,s=2,p=1)+ Relu() \\ 
        ConvTranspose2d(i=4,o=4,k=3,s=2,p=1)+Relu()  \\ 
        ConvTranspose2d(i=4,o=4,k=3,s=2,p=1)+Relu()  \\ 
        ConvTranspose2d(i=4,o=4,k=3,s=2,p=1)+Relu()  \\ 
        ConvTranspose2d(i=4,o=4,k=3,s=2,p=1)+Relu()  \\ 
        ConvTranspose2d(i=4,o=4,k=3,s=2,p=1)+Relu()  \\ 
        ConvTranspose2d(i=4,o=2,k=3,s=2,p=2)         \\ 
        AdaptiveAvgPool2d($n$) \\
        \bottomrule
        \end{tabular}
\end{table}

\subsection{Convection-diffusion equation}\label{appen:ad}
\begin{table}[h]
    \centering
    \caption{Training hyperparameters for convection-diffusion equation.}
    \label{tab:ad_hyper}
    \begin{tabular}{@{}lccccc@{}}
        \toprule
        & learning rate & batch size & $K$         & xavier init & grad clip \\ \midrule
FNS(Jacobi) & $10^{-4}$     & $100$      & $10$        & $10^{-2}$   & false     \\
FNS(Conv)   & $10^{-4}$     & $100$      & $1\sim 100$ & $10^{-2}$   & 1.0      \\ \bottomrule
\end{tabular}
\end{table}
\begin{table}[h]
    \centering
    \caption{HyperNN architecture parameters for convection-diffusion equation. Notations in ConvTranspose2d are: i: in channels; o: out channels; k: kernel size; s: stride; p: padding.}
    \label{tab:ad_hyper_nn}
        \begin{tabular}{@{}cl@{}}
        \toprule
        $\Phi$                       & Hyper NN                                     \\ \midrule
        Conv(i=1,o=8, k=3, s=2, p=1) & ConvTranspose2d(i=1,o=4,k=3,s=2,p=1)+Relu() \\ 
        Conv(i=8,o=1, k=3, s=2, p=1) & ConvTranspose2d(i=4,o=4,k=3,s=2,p=1)+Relu()  \\ 
                                     & ConvTranspose2d(i=4,o=4,k=3,s=2,p=1)+Relu()  \\ 
                                     & ConvTranspose2d(i=4,o=4,k=3,s=2,p=1)+Relu()  \\ 
                                     & ConvTranspose2d(i=4,o=2,k=3,s=2,p=2)         \\ \bottomrule
        \end{tabular}
\end{table}

\subsection{Helmholtz equation}\label{appen:helm}
\begin{table}[h]
    \centering
    \caption{Training hyperparameters for Helmholtz equation.}
    \label{tab:helm_hyper}
    \begin{tabular}{@{}lccccc@{}}
        \toprule
        & learning rate & batch size & $K$         & xavier init & grad clip \\ \midrule
FNS(Krylov)    & $10^{-4}$     & $100$      & $1\sim 100$ & $10^{-2}$   & 1.0      \\ \bottomrule
\end{tabular}
\end{table}
\begin{table}[h]
    \centering
    \caption{HyperNN architecture parameters for Helmholtz equation. Notations in ConvTranspose2d are: i: in channels; o: out channels; k: kernel size; s: stride; p: padding.}
    \label{tab:helm_hyper_nn}
    \begin{tabular}{@{}l@{}}
        \toprule
        ConvTranspose2d(i=1,o=4,k=3,s=2,p=1)+ Relu() \\ 
        ConvTranspose2d(i=4,o=4,k=3,s=2,p=1)+Relu()  \\ 
        ConvTranspose2d(i=4,o=4,k=3,s=2,p=1)+Relu()  \\ 
        ConvTranspose2d(i=4,o=4,k=3,s=2,p=1)+Relu()  \\ 
        ConvTranspose2d(i=4,o=4,k=3,s=2,p=1)+Relu()  \\ 
        ConvTranspose2d(i=4,o=4,k=3,s=2,p=1)+Relu()  \\ 
        ConvTranspose2d(i=4,o=2,k=3,s=2,p=2)         \\ 
        AdaptiveAvgPool2d($n$) \\
        \bottomrule
        \end{tabular}
\end{table}

\bibliographystyle{unsrt}
\bibliography{fns}

\begin{thebibliography}{10}

\bibitem{barrett1994templates}
Richard Barrett, Michael Berry, Tony~F Chan, James Demmel, June Donato, Jack
  Dongarra, Victor Eijkhout, Roldan Pozo, Charles Romine, and Henk Van~der
  Vorst.
\newblock {\em Templates for the solution of linear systems: building blocks
  for iterative methods}.
\newblock SIAM, 1994.

\bibitem{saad2003iterative}
Yousef Saad.
\newblock {\em Iterative methods for sparse linear systems}.
\newblock SIAM, 2003.

\bibitem{hestenes1952methods}
Magnus~R Hestenes and Eduard Stiefel.
\newblock Methods of conjugate gradients for solving.
\newblock {\em Journal of research of the National Bureau of Standards},
  49(6):409, 1952.

\bibitem{saad1986gmres}
Youcef Saad and Martin~H Schultz.
\newblock Gmres: A generalized minimal residual algorithm for solving
  nonsymmetric linear systems.
\newblock {\em SIAM Journal on scientific and statistical computing},
  7(3):856--869, 1986.

\bibitem{brandt1977multi}
Achi Brandt.
\newblock Multi-level adaptive solutions to boundary-value problems.
\newblock {\em Mathematics of computation}, 31(138):333--390, 1977.

\bibitem{briggs2000multigrid}
William~L Briggs, Van~Emden Henson, and Steve~F McCormick.
\newblock {\em A multigrid tutorial}.
\newblock SIAM, 2000.

\bibitem{trottenberg2000multigrid}
Ulrich Trottenberg, Cornelius~W Oosterlee, and Anton Schuller.
\newblock {\em Multigrid}.
\newblock Elsevier, 2000.

\bibitem{falgout2006introduction}
Robert~D Falgout.
\newblock An introduction to algebraic multigrid.
\newblock Technical report, Lawrence Livermore National Lab.(LLNL), Livermore,
  CA (United States), 2006.

\bibitem{xu2017algebraic}
Jinchao Xu and Ludmil Zikatanov.
\newblock Algebraic multigrid methods.
\newblock {\em Acta Numerica}, 26:591--721, 2017.

\bibitem{hsieh2019learning}
Jun-Ting Hsieh, Shengjia Zhao, Stephan Eismann, Lucia Mirabella, and Stefano
  Ermon.
\newblock Learning neural pde solvers with convergence guarantees.
\newblock {\em arXiv preprint arXiv:1906.01200}, 2019.

\bibitem{luna2021accelerating}
Kevin Luna, Katherine Klymko, and Johannes~P Blaschke.
\newblock Accelerating gmres with deep learning in real-time.
\newblock {\em arXiv preprint arXiv:2103.10975}, 2021.

\bibitem{weymouth2021data}
Gabriel~D Weymouth.
\newblock Data-driven multi-grid solver for accelerated pressure projection.
\newblock {\em arXiv preprint arXiv:2110.11029}, 2021.

\bibitem{tomasi2021construction}
Claudio Tomasi and Rolf Krause.
\newblock Construction of grid operators for multilevel solvers: a neural
  network approach.
\newblock {\em arXiv preprint arXiv:2109.05873}, 2021.

\bibitem{taghibakhshi2021optimization}
Ali Taghibakhshi, Scott MacLachlan, Luke Olson, and Matthew West.
\newblock Optimization-based algebraic multigrid coarsening using reinforcement
  learning.
\newblock {\em Advances in Neural Information Processing Systems},
  34:12129--12140, 2021.

\bibitem{huang2021learning}
Ru~Huang, Ruipeng Li, and Yuanzhe Xi.
\newblock Learning optimal multigrid smoothers via neural networks.
\newblock {\em arXiv preprint arXiv:2102.12071}, 2021.

\bibitem{wang4110904learning}
Fan Wang, Xiang Gu, Jian Sun, and Zongben Xu.
\newblock Learning-based local weighted least squares for algebraic multigrid
  method.
\newblock {\em Available at SSRN 4110904}.

\bibitem{fanaskov2021neural}
Vladimir Fanaskov.
\newblock Neural multigrid architectures.
\newblock In {\em 2021 International Joint Conference on Neural Networks
  (IJCNN)}, pages 1--8. IEEE, 2021.

\bibitem{chen2022meta}
Yuyan Chen, Bin Dong, and Jinchao Xu.
\newblock Meta-mgnet: Meta multigrid networks for solving parameterized partial
  differential equations.
\newblock {\em Journal of Computational Physics}, 455:110996, 2022.

\bibitem{katrutsa2020black}
Alexandr Katrutsa, Talgat Daulbaev, and Ivan Oseledets.
\newblock Black-box learning of multigrid parameters.
\newblock {\em Journal of Computational and Applied Mathematics}, 368:112524,
  2020.

\bibitem{greenfeld2019learning}
Daniel Greenfeld, Meirav Galun, Ronen Basri, Irad Yavneh, and Ron Kimmel.
\newblock Learning to optimize multigrid pde solvers.
\newblock In {\em International Conference on Machine Learning}, pages
  2415--2423. PMLR, 2019.

\bibitem{luz2020learning}
Ilay Luz, Meirav Galun, Haggai Maron, Ronen Basri, and Irad Yavneh.
\newblock Learning algebraic multigrid using graph neural networks.
\newblock In {\em International Conference on Machine Learning}, pages
  6489--6499. PMLR, 2020.

\bibitem{antonietti2021accelerating}
Paola~F Antonietti, Matteo Caldana, and Luca Dede.
\newblock Accelerating algebraic multigrid methods via artificial neural
  networks.
\newblock {\em arXiv preprint arXiv:2111.01629}, 2021.

\bibitem{stanziola2021helmholtz}
Antonio Stanziola, Simon~R Arridge, Ben~T Cox, and Bradley~E Treeby.
\newblock A helmholtz equation solver using unsupervised learning: Application
  to transcranial ultrasound.
\newblock {\em Journal of Computational Physics}, 441:110430, 2021.

\bibitem{kapturowski2018recurrent}
Steven Kapturowski, Georg Ostrovski, John Quan, Remi Munos, and Will Dabney.
\newblock Recurrent experience replay in distributed reinforcement learning.
\newblock In {\em International conference on learning representations}, 2018.

\bibitem{azulay2022multigrid}
Yael Azulay and Eran Treister.
\newblock Multigrid-augmented deep learning preconditioners for the helmholtz
  equation.
\newblock {\em arXiv preprint arXiv:2203.11025}, 2022.

\bibitem{erlangga2006novel}
Yogi~A Erlangga, Cornelis~W Oosterlee, and Cornelis Vuik.
\newblock A novel multigrid based preconditioner for heterogeneous helmholtz
  problems.
\newblock {\em SIAM Journal on Scientific Computing}, 27(4):1471--1492, 2006.

\bibitem{calandra2012flexible}
Henri Calandra, Serge Gratton, Julien Langou, Xavier Pinel, and Xavier Vasseur.
\newblock Flexible variants of block restarted gmres methods with application
  to geophysics.
\newblock {\em SIAM Journal on Scientific Computing}, 34(2):A714--A736, 2012.

\bibitem{um2020solver}
Kiwon Um, Robert Brand, Yun~Raymond Fei, Philipp Holl, and Nils Thuerey.
\newblock Solver-in-the-loop: Learning from differentiable physics to interact
  with iterative pde-solvers.
\newblock {\em Advances in Neural Information Processing Systems},
  33:6111--6122, 2020.

\bibitem{nikolopoulos2022ai}
Stefanos Nikolopoulos, Ioannis Kalogeris, Vissarion Papadopoulos, and George
  Stavroulakis.
\newblock Ai-enhanced iterative solvers for accelerating the solution of large
  scale parametrized linear systems of equations.
\newblock {\em arXiv preprint arXiv:2207.02543}, 2022.

\bibitem{stanaityte2020ilu}
Rita Stanaityte et~al.
\newblock {\em ILU and Machine Learning Based Preconditioning For The
  Discretized Incompressible Navier-Stokes Equations.}
\newblock PhD thesis, 2020.

\bibitem{kaneda2022deep}
Ayano Kaneda, Osman Akar, Jingyu Chen, Victoria Kala, David Hyde, and Joseph
  Teran.
\newblock A deep gradient correction method for iteratively solving linear
  systems.
\newblock {\em arXiv preprint arXiv:2205.10763}, 2022.

\bibitem{margenberg2022neural}
Nils Margenberg, Dirk Hartmann, Christian Lessig, and Thomas Richter.
\newblock A neural network multigrid solver for the navier-stokes equations.
\newblock {\em Journal of Computational Physics}, 460:110983, 2022.

\bibitem{margenberg2021deep}
Nils Margenberg, Robert Jendersie, Thomas Richter, and Christian Lessig.
\newblock Deep neural networks for geometric multigrid methods.
\newblock {\em arXiv preprint arXiv:2106.07687}, 2021.

\bibitem{cooley1965algorithm}
James~W Cooley and John~W Tukey.
\newblock An algorithm for the machine calculation of complex fourier series.
\newblock {\em Mathematics of computation}, 19(90):297--301, 1965.

\bibitem{sleijpen1993bicgstab}
Gerard~LG Sleijpen and Diederik~R Fokkema.
\newblock Bicgstab (ell) for linear equations involving unsymmetric matrices
  with complex spectrum.
\newblock {\em Electronic Transactions on Numerical Analysis.}, 1:11--32, 1993.

\bibitem{swarztrauber1977methods}
Paul~N Swarztrauber.
\newblock The methods of cyclic reduction, fourier analysis and the facr
  algorithm for the discrete solution of poisson’s equation on a rectangle.
\newblock {\em Siam Review}, 19(3):490--501, 1977.

\bibitem{kingma2014adam}
Diederik~P Kingma and Jimmy Ba.
\newblock Adam: A method for stochastic optimization.
\newblock {\em arXiv preprint arXiv:1412.6980}, 2014.

\bibitem{paszke2019pytorch}
Adam Paszke, Sam Gross, Francisco Massa, Adam Lerer, James Bradbury, Gregory
  Chanan, Trevor Killeen, Zeming Lin, Natalia Gimelshein, Luca Antiga, et~al.
\newblock Pytorch: An imperative style, high-performance deep learning library.
\newblock {\em Advances in neural information processing systems}, 32, 2019.

\bibitem{huang2017densely}
Gao Huang, Zhuang Liu, Laurens Van Der~Maaten, and Kilian~Q Weinberger.
\newblock Densely connected convolutional networks.
\newblock In {\em Proceedings of the IEEE conference on computer vision and
  pattern recognition}, pages 4700--4708, 2017.

\bibitem{golub1961chebyshev}
Gene~H Golub and Richard~S Varga.
\newblock Chebyshev semi-iterative methods, successive overrelaxation iterative
  methods, and second order richardson iterative methods.
\newblock {\em Numerische Mathematik}, 3(1):157--168, 1961.

\bibitem{adams2003parallel}
Mark Adams, Marian Brezina, Jonathan Hu, and Ray Tuminaro.
\newblock Parallel multigrid smoothing: polynomial versus gauss--seidel.
\newblock {\em Journal of Computational Physics}, 188(2):593--610, 2003.

\bibitem{mises1929praktische}
RV~Mises and Hilda Pollaczek-Geiringer.
\newblock Praktische verfahren der gleichungsaufl{\"o}sung.
\newblock {\em ZAMM-Journal of Applied Mathematics and Mechanics/Zeitschrift
  f{\"u}r Angewandte Mathematik und Mechanik}, 9(1):58--77, 1929.

\end{thebibliography}

\end{document}